\documentclass{article}

\usepackage{parskip}

\usepackage{color}
\usepackage{graphicx}
\usepackage{multirow,multicol}
\usepackage{amssymb,amsmath,amsthm,amsfonts}
\usepackage{tikz}
\usepackage{tikz-3dplot}

\theoremstyle{definition}
\newtheorem{example}{Example}[section]

\title{Uncertainty quantification of coal seam gas production prediction using Polynomial Chaos}
\author{Thomas A. McCourt\footnote{This research has been conducted with funding and support from the University of Queensland Centre for Coal Seam Gas and its industry members (APLNG, Arrow Energy, QGC (Shell) and Santos).
We acknowledge CMG Ltd. for donating commercial reservoir simulation software. Foundation CMG is acknowledged for partial funding of this study.} \and Suzanne Hurter\footnotemark[1] \and Brodie A. J. Lawson\footnotemark[1] \and Fengde Zhou\footnotemark[1] \and Bevan Thompson\footnotemark[1] \and Stephen Tyson\footnotemark[1] \and Diane Donovan\footnotemark[1]}

\newcommand{\cM}{\mathcal{M}}
\newcommand{\cP}{\mathcal{P}}
\newcommand{\cL}{\mathcal{K}}
\newcommand{\cN}{\mathcal{N}}
\newcommand{\cI}{\mathcal{I}}

\begin{document}

\maketitle

\begin{abstract}
A surrogate model approximates a computationally expensive solver. Polynomial Chaos is a method used to construct surrogate models by summing combinations of carefully chosen polynomials. The polynomials are chosen to respect the probability distributions of the uncertain input variables (parameters); this allows for both uncertainty quantification and global sensitivity analysis. 

In this paper we apply these techniques to a commercial solver for the estimation of peak gas rate and cumulative gas extraction from a coal seam gas well. The polynomial expansion is shown to honour the underlying geophysics with low error when compared to a much more complex and computationally slower commercial solver. We make use of advanced numerical integration techniques to achieve this accuracy using relatively small amounts of training data.
\end{abstract}

\begin{center}
\begin{tabular}{l}
and Global Sensitivity Analysis
\end{tabular}
\end{center}

\section{Introduction}
In coal seam gas production reservoir simulations play an important role in estimating extraction rates and associated economic forecasting \cite{Aminian2009, Scott2008, Zhou2014}.
  Under appropriate pressure and temperature conditions  buried peat, gas is generated by thermal or biogenic processes, and a network of fractures (face and butt cleats) is formed \cite{Gayer1996}. %In well developed cleat systems, this fracture network presents as a network of face and butt cleats.

    Prior to production, most of the coal seam gas (CSG, predominantly methane) exists as adsorbate in the micro-pores \cite{Gray1987} while cleats are fully or partially saturated with water.  Dewatering is used to decrease the pressure in the cleats and when this pressure is less than the critical desorption pressure the gas desorbs from the matrix into the cleats \cite{Seidle2011}. Once the gas saturation exceeds the residual gas saturation, gas begins to flow along with water to a producer. This  process is assumed to obey Darcy's Law \cite{Aminian2009,Zhou2014}  and is simulated numerically  to predict gas and water production for CSG wells.  The gas and water production increases until the production rates reach a peak (not simultaneous for water and gas) and decline thereafter. Over time a  reduction in the gas content and pressure in the matrix results in a decline in gas production rates.

  Key indicators are used to predict the gas production curve, for instance: peak gas production rate, the time to peak gas rate, and other production decline coefficients \cite{Aminian2009,Keim2011}.    Aminian and Ameri \cite{Aminian2009} derived a formula for calculating the peak gas rate based on the fluid radial flow equation.  Bhavsar \cite{Bhavsar2005} estimated the dimensionless peak gas rate using a regression equation with critical gas desorption pressure, skin factor, porosity, Langmuir volume and Langmuir pressure.  Zhou \cite{Zhou2014} developed equations to predict the peak gas rate, peak gas rate arrival time and decline rate by multiple regression of 200 simulations.  Those simulations were based on a static model with varied skin factor, porosity, permeability, model geometry, desorption time, Langmuir volume and pressure, thickness, dewatering pressure, and critical desorption pressure.

For the operation of a CSG field, %timing of 
peak gas and the cumulative gas production over a period of time are important measures of field performance and aid in field management decisions. Uncertainty in these model results has significant impact on economic and budgeting considerations.

Sophisticated simulation packages have been developed to capture the complicated nature of the pressure changes and  fluid flow regimes in hydrocarbon reservoirs. Techniques such as experimental design and Monte Carlo simulations have been incorporated into some of the commercial simulation packages and applied across various industries; see for instance \cite{Collins2015,Yeten2005} for studies on petroleum recovery.   Li, Sarma and Zhang \cite{Li2011} used  probabilistic collocation methods to quantify uncertainty in petroleum reservoirs  and then compared these techniques with the more traditional experimental design approach. A key disadvantage of the experimental design approach is that the stochastic nature of the problem may not be captured as the probability distributions of the parameters are ignored. While Monte Carlo simulations can overcome this problem, they can be time consuming and expensive, especially for production optimisation and uncertainty quantification. In addition,  the complexity of these simulation packages renders them virtual `black-boxes' in the sense that modifications to their inner workings are not possible.
 In this setting techniques that reduce computational cost, while at the same time deliver  additional statistical information, are invaluable.

As an alternative
Sarma and Xie \cite{Sarma2011} proposed  the use  of Polynomial Chaos Expansions (PCEs) to develop simulations for forecasting oil reservoir performance. Sarma and Xie treated the reservoir simulation as a black-box, approximating it with a PCE.

Polynomial Chaos Expansion (PCE) is a mathematical technique for taking very complex models  and constructing \textit{surrogate} models that take the same inputs and accurately and  more efficiently approximate the outputs. \textit{Non-intrusive polynomial chaos} is a variation of this technique that requires no modification of the original model, but can replace a  virtual `black-box' with a simple surrogate that behaves similarly in terms of its inputs and outputs. This much simpler model can then be evaluated very quickly, allowing for better exploration of the parameter space, uncertainty quantification, and, potentially, better forecasting tools and new work-flows that take advantage of the significant reduction in computational time. Moreover, it enables the ready extraction of Sobol' indices and hence provides a global sensitivity analysis for the surrogate, which in turn approximates the equivalent analysis for the original model.

Further examples of the application of PCEs in the study of subsurface oil reservoirs can be found in \cite{Alkhatib2014} and \cite{Jansen2008} where PCEs are used to study surfactant-polymer flooding for oil recovery processes. %{\color{red}
The application PCEs are also applied to subsurface oil reservoir models in \cite{Bazargan2013} and \cite{Dai2014}; indeed \cite{Dai2014} uses the constructed PCEs to perform global sensitivity analysis in a similar manner to that used later in this paper for coal seam gas models. 
The literature also documents  the use of PCEs  in many other applications, for instance, to study the trapping of CO$_2$ \cite{Babaei2015b,Oladyshkin2011}, flow in porous media \cite{Fajraoui2011}, and subsurface flows \cite{Babaei2015a,Elsheikh2014}. %{\color{red} 
Applications of PCE techniques to CSG modelling have recently begun to be addressed in the literature \cite{Senthamaraikkannan2016}, but are not used in the conventional CSG industry, supporting the argument of a need to study and apply PCE to problems unique to CSG.

In the current paper we wish to take advantage of the power of PCE and use it to emulate a commercial simulation package for the estimation of  peak and total gas production for a single well tapping water and gas from a single coal seam.
 The production forecasts for coal seam gas are heavily dependent on the time of peak gas arrival, a problem somewhat  unique to coal seam gas. The two phase flow of gas and water behaves differently from two phase flow of oil and water, with the peak rates occurring at different times for water and gas with the gas flowing after the water is much diminished, as demonstrated in Figure \ref{fig:peak-gas-water}. As evidenced by current industry practices, there is a need to investigate models developed specifically for CSG; a need that indicates the timeliness of and motivation for this current study. 

Here the commercial  package is treated as a black-box and is used to generate training points, on which the surrogate models are based. Further points are chosen via Latin Hypercube sampling across the parameter space, and  used to provide error estimates comparing  the surrogate and the commercial solver. In addition, we extract the Sobol' indices and discuss the sensitivity of the surrogate models to both individual and pair-wise model parameters.
The objective is to investigate how well PCE techniques emulate the workings of the Compute Modelling Group GEM  tool \cite{CMG2014}.  We have made a conscious decision to  model a single well in a multi-cell field with homogeneous initial values. To introduce a heterogeneous multi-cell field would introduce tens of thousand of parameters and require significant dimensionality reduction \cite{Laloy2013,Marzouk2009,Zhang2004} before constructing a PCE surrogate model. This added complexity leads to the risk of complicating the analysis of the accuracy of the PCE surrogate model, as well as the effectiveness of the derived Sobol' indices for sensitivity analysis.
Thus in this scoping study we have `up-scaled'  
the parameter values for permeability, porosity and Langmuir parameters while maintaining a fine grid structure, 
so we may accurately provide an assessment of the PCE technique.
A discussion of the choice of uncertain parameters and their distributions is left to Section \ref{sec:model} where we provide background on the commercial simulation package and the single well model.

 In Section \ref{sec:pce} we detail the theory of spectral expansion using Polynomial Chaos and the related numerical integration using full grid Gauss-Legendre quadrature and sparse grid Clenshaw-Curtis quadrature.

 Section \ref{subsubsec:GSA} provides details on the theory of Sobol' indices for sensitivity analysis and their direct calculation from PCEs. As Sudret \cite{Sudret2008} states, Sobol' indices do not assume linearity or monotonicity in the model, therefore they  are good descriptors of the sensitivity of the model to its input parameters. To calculate the full description of parameter interaction using Monte Carlo simulations  would require many evaluations, which quickly becomes infeasible. Thus the derivation of the Sobol' indices directly from the PCE can provide significant computational savings \cite{Sudret2008}.

 Section \ref{sec:pce_here} describes how these techniques have been applied to the commercial solver and includes error estimates, uncertainty quantification and global sensitivity analysis. We highlight that the PCEs allow for such analysis to be performed in a computationally inexpensive manner.

\section{The Model}
\label{sec:model}

In this paper we will demonstrate the application of PCE as an emulator for the commercial simulation package  GEM, developed by the Computer Modelling Group (CMG)  \cite{CMG2014}, that is used to study of  the extraction of adsorbed gases.

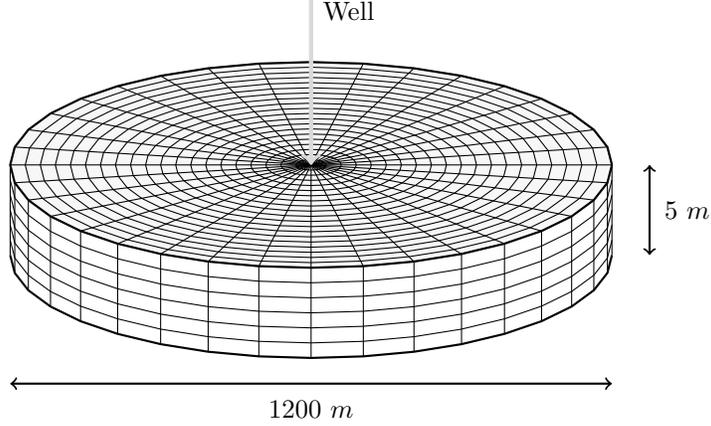
\begin{figure}
\begin{center}
%\resizebox{\linewidth}{!}{
\tdplotsetmaincoords{70}{0}
\begin{tikzpicture}[tdplot_main_coords]
\def\RI{4}
\def\RII{4}

\draw[thick,<->] (4.5,0) -- (4.5, 1.2cm );
\node at (5,0.6cm)   { $5\ m$};

\draw[thick,<->] (-4,-5) -- (4, -5);
\node at (0,-6)   { $1200\ m$};

\draw[thick] (\RI,0)
  \foreach \x in {0, 350, 340, 330, 320, 310, 300, 290, 280, 270, 260, 250, 240, 230, 220, 210, 200, 190, 180} { --  (\x:\RI) node at (\x:\RI) (R1-\x) {} };

\begin{scope}[yshift=1.2cm]
\draw[thick,fill=gray!30,opacity=0.2] (\RII,0)
  \foreach \x in {0, 10, 20, 30, 40, 50, 60, 70, 80, 90, 100, 110, 120, 130, 140, 150, 160, 170, 180, 190, 200, 210, 220, 230, 240, 250, 260, 270, 280, 290, 300, 310, 320, 330, 340, 350, 360} { --  (\x:\RII) node at (\x:\RII) (R2-\x) {}};
\end{scope}

%%%%%%%%%%%%%%%%%%%%%%%%%%%%%%%

% Vertical lines
\foreach \x in {0, 350, 340, 330, 320, 310, 300, 290, 280, 270, 260, 250, 240, 230, 220, 210,200, 190, 180} { \draw (R1-\x.center) -- (R2-\x.center); };

\begin{scope}[yshift=1.2cm]
\foreach \x in {0, 350, 340, 330, 320, 310, 300, 290, 280, 270, 260, 250, 240, 230, 220, 210,200, 190, 180, 170, 160, 150, 140, 130, 120, 110, 100, 90, 80, 70, 60, 50, 40, 30, 20, 10} { \draw (R2-\x.center) -- (0,0); };
\end{scope}

%%%%%%%%%%%%%%%%%%%%%%%%%%%%%%%

%Layers
\foreach \y in {0.2, 0.4, 0.6, 0.8, 1.0}{
\begin{scope}[yshift=\y cm]
\draw (\RI,0)
  \foreach \x in {0, 350, 340, 330, 320, 310, 300, 290, 280, 270, 260, 250, 240, 230, 220, 210, 200, 190, 180} { --  (\x:\RI) node at (\x:\RI) (R1-\x) {} };
\end{scope}
}

%%%%%%%%%%%%%%%%%%%%%%%%%%%%%%%

%concentric stuff

\foreach \r in {0.2, 0.4, 0.6, 0.8, 1, 1.2, 1.4, 1.6, 1.8, 2, 2.2, 2.4, 2.6, 2.8, 3, 3.2, 3.4, 3.6, 3.8}
{
\begin{scope}[yshift=1.2 cm]
\draw (\r,0)
  \foreach \x in {360, 350, 340, 330, 320, 310, 300, 290, 280, 270, 260, 250, 240, 230, 220, 210, 200, 190, 180, 170, 160, 150, 140, 130, 120, 110, 100, 90, 80, 70, 60, 50, 40, 30, 20, 10, 0} { --  (\x:\r) node at (\x:\r) (R1-\x) {} };
\end{scope}
}

\begin{scope}[yshift=1.2 cm]
\draw[thick] (\RI,0)
  \foreach \x in {360, 350, 340, 330, 320, 310, 300, 290, 280, 270, 260, 250, 240, 230, 220, 210, 200, 190, 180, 170, 160, 150, 140, 130, 120, 110, 100, 90, 80, 70, 60, 50, 40, 30, 20, 10, 0} { --  (\x:\RI) node at (\x:\RI) (R1-\x) {} };
\end{scope}

%%%%%%%%%%%%%%%%%%%%%%%%%%%%%%%

\draw[ultra thick, color=gray!30, ->] (0,10) -- (0,1.2cm);
\node at (0.5,9.5) {Well};
\end{tikzpicture}
%}
\end{center}
\caption{Illustration of the grid system with radial distance plotted on a logarithmic scale.}
\label{fig:a_cell}
\end{figure}

A single well in a homogeneous field is simulated using the static model developed by Zhou in \cite{Zhou2014} and is then used to predict peak and total gas rates. In this setting the model geometry is  a fixed radial grid system, referenced by radial distance $x$, bearing (angle) $y$ and depth (thickness) $z$, with the radial distance  equally partitioned into 40 segments on a logarithmic scale, the bearing partitioned into 36 segments and the depth into 6 layers, as illustrated in Figure \ref{fig:a_cell}.
The radial size is $600\ m$ and is of the order of inter-well distances of vertical wells in the Surat Basin, eastern Australia \cite{Lagendijk2010}.  
The coal seam thickness is $5\ m$ which is about the average coal seam thickness in the Upper Juandah of the Surat Basin, eastern Australia \cite{Martin2013}. One well is located at the centre and perforates  all six layers. The top depth of this model is $440\ m$ which is similar to the average burial depth of the coal of the Upper Juandah.
In this set-up, the field is modelled by 8640 grid cells of increasing volume as we move away from the central well. The flow is modelled  cell by cell with flow between neighbouring cells individually accounted for in the model.

A  two phase flow involving water and gas is assumed, with the gas  not dissolved in the water.	The temperature is constant and the coal cleat distribution is isotropic.
The gas desorption obeys the Langmuir equation.
Given these conditions, the basic mass conservation equation  for gas is given in Equation \eqref{one} and for water by Equation \eqref{two}.
\begin{eqnarray}
\nabla
\cdot\left[\frac{\rho_g kk_{g}}{\mu_g}\left(\nabla p_{gf}-9.81\rho_g \nabla H\right)\right]-q_{gf}+q_{gmf}&=&\frac{\partial}{\partial t}\left(\phi_f\rho_g S_{gf}\right),\label{one}\\
\nabla
\cdot\left[\frac{\rho_w kk_{w}}{\mu_g}\left(\nabla p_{wf}- 9.81\rho_w\nabla H\right)\right]-q_{wf}&=&\frac{\partial}{\partial t}\left(\phi_f\rho_w S_{wf}\right),\label{two}
%P_{cgwf}\left(S_{gf}\right)=p_{gf}-p_{wf}&\mbox{and}&S_{gf}+S_{wf}=1,
\end{eqnarray}
where $S_{gf}+S_{wf}=1$, subscripts $w,\ g,\ f$ and $m$ refer to water, gas, fracture, and matrix, respectively and $0$ indicating an initial condition.   Other parameters are: the density, $\rho$; the viscosity, $\mu$; the flow rate, $q$; the porosity, $\phi$; the saturation, $S$; the pressure, $p$; the elevation, $H$; the time, $t$; and the capillary pressure, $P_c$.
The absolute permeability, $k$,  is calculated as in \cite{McKee1987} by
\begin{eqnarray}
\frac{k}{k_0}=\left(\frac{\phi}{\phi_0}\right)^a\label{thirteen}
 \end{eqnarray}
 where  the exponent $a$ is set to $3$ for naturally-fractured reservoirs like coals.
The gas diffusion rate, $q_{gmf}$, between the cleats and the matrix is calculated by
\begin{eqnarray}
q_{gmf}&=&\frac{\rho_g p_{av}V_m}{\tau}F(S_g)\left(\frac{V_Lp_{mo}}{P_L+p_{mo}}-\frac{V_Lp_m}{P_L+p_m}\right),
\end{eqnarray}
where, $P_L$ and $V_L$ are the Langmuir pressure and Langmuir volume, respectively; $\tau$ is the desorption time of gas; and $V_m$ is the coal volume of a cell. An average pressure between the initial matrix pressure and the actual matrix pressure is $p_{av} = (p_{m0}+p_m)/2$.  The quantity $F(S_g)$ is a function of gas saturation in the fractures and models water blocking damage.  The desorption time, $\tau$ equals $1/(\sigma \cdot D)$,
where $\sigma$ is the shape factor and $D$ is the diffusion coefficient in $m^2$ per day. The  Gilman and Kazemi equation, $\sigma=4\cdot \left(1/L_x^2+1/L_y^2+1/L_z^2\right)$, is used to calculate the shape factor, where $L$ is the fracture (cleat) spacing in 
the radial grid indices $x$, $y$ and $z$.

\begin{figure}
\begin{center}

\resizebox{\linewidth}{!}{\begin{tabular}{c}
\includegraphics[trim=0mm 30mm 95mm 0mm,clip,width = 14cm]{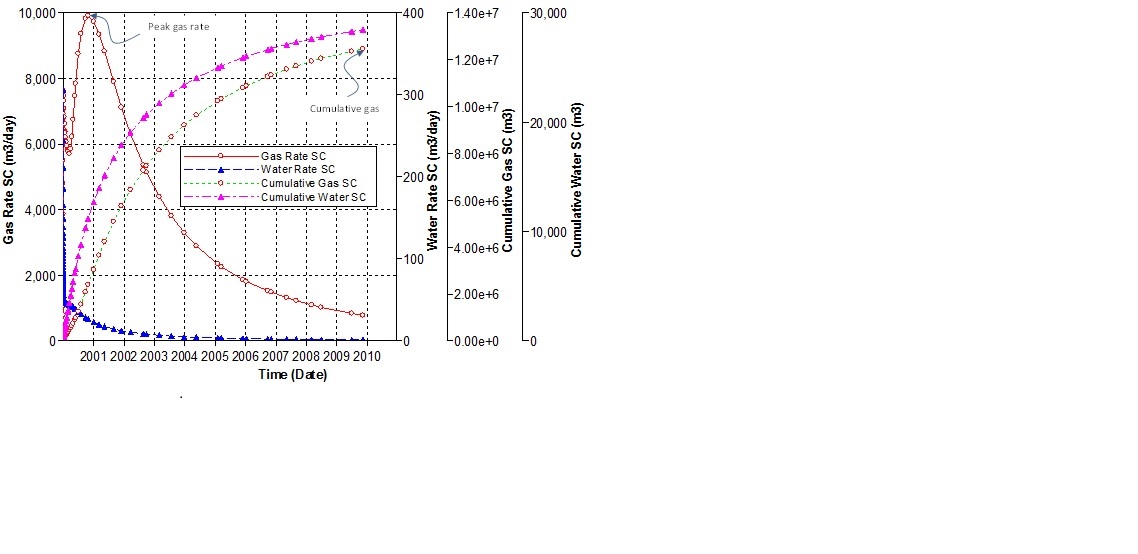}
\end{tabular}}

\end{center}
\caption{Simulated gas and water production rates and cumulative gas and water production for one run of the model. The permeability, porosity, Langmuir volume and pressure for this case are 43 mD, 0.7\%, 0.255 mol/kg, and 5882 kPa, respectively. This case was simulated with well bottom-hole pressure at 300kPa as the constraint.}
\label{fig:peak-gas-water}
\end{figure}

Industry experience, in the Surat Basin, indicates that the most significant effect on field production performance is  variability and uncertainty in the permeability and porosity. Porosity is related to the amount of water produced, affecting the required dimensionality of water treatment plants, while permeability affects the total flow of water and gas. Langmuir volume ($V_L$) is the maximum volume of gas a coal can adsorb onto its surface while  Langmuir pressure ($P_L$) is the pressure at which the storage capacity of a coal is equal to half the Langmuir volume. Langmuir volume relates to the maximum gas stored for a specified volume of coal (i.e. gas-in-place) while Langmuir pressure controls the rate at which gas can be recovered per unit pressure drop. A high $P_L$ implies faster recovery of gas per unit pressure drop and ultimately higher cumulative gas production \cite{Bahrami2015}. These two parameters are critical in determining the arrival time of the peak gas rate.  The time to peak gas arrival impacts significantly on production forecasts for coal seam gas, with the two phase flow of gas and water showing the peak rates occur at different times, as demonstrated in Figure \ref{fig:peak-gas-water}.

Given that this scoping study focusses on a static model in a homogeneous single well field the four parameters permeability, porosity, Langmuir volume and Langmuir pressure are deemed the critical parameters for the determination of peak gas and so the PCE will be built to explore uncertainty and sensitivity of the model to these parameters.  The values of the fracture porosity, fracture permeability, Langmuir volume and Langmuir pressure are uncertain and we refer to them as the \textit{input variables} of the model.

Although other parameters influence the production behaviour of coal (e.g. relative permeability, heterogeneity, stress, coal composition, amongst others), Langmuir parameters, permeability and porosity are prioritised and selected for this application based on the impact that industry attributes to these for economic valuation. 

 The remaining model properties will be assigned fixed values as listed in Table \ref{tab:inputs}. The initial pressure in the cleats is $4440\ kPa$ at a reference depth of $440\ m$  assuming hydrostatic pressure, while the initial pressure in the matrix is $2750\ kPa$. This leads to an initial gas saturation of about $77\%$  in the matrix. The above values, properties and ranges listed in Table \ref{tab:inputs} and a desorption time of 0.4 days are sourced from \cite{Scott2008} and communications with an Australian CSG company. The boundary conditions for the models are taken as closed for top, bottom and lateral boundaries.

\begin{table}
\begin{center}
\resizebox{\linewidth}{!}{
\begin{tabular}{lll}
\hline
\textbf{Parameters}	&\textbf{Ranges}&	\textbf{Units}\\
\hline
Model type &Dual porosity	&-\\
Shape factor formulation	&Gilman-Kazemi	&-\\
Fluid component model	&Peng-Robinson equation	&-\\
Model geometry	&Radial grid	&-\\
Grid system	&$40\times 36\times 6$	&-\\
Thickness $h$	&$5$&$m$\\
Matrix porosity $\phi _m$	&$0.01$ 	&$1$\\
Fracture porosity $\phi$&	$[0.005,0.05]$	&$1$\\
Matrix permeability $k_m$&	$0.01$&	$mD$\\
Fracture permeability $k_x$&	$[10,1000]$ &	$mD$\\
Initial reservoir pressure in fracture $P_{rf}$ &	$4440$&	$kPa$\\
Initial reservoir pressure in matrix $P_{rm}$&	$2750$ &	$kPa$\\
Reservoir temperature $T_r$ &	$45$ &	$\,^o C$\\
Reciprocal of the Langmuir pressure $b$ &	$[0.00017,0.0003]$ &	$1/kPa$\\
Langmuir volume $V_L$ &	$[0.2,1]$ &	$gmole/kg$\\
Sorption time $\tau$&	$0.4$&	days\\
Coal density $\rho$	&$1435$&	$kg/m^3$\\
Well bottom-hole pressure & $101.3$ & $kPa$\\
\hline	
\end{tabular}}
\caption{Reservoir properties and ranges used in simulation.}\label{tab:inputs}
\end{center}
\end{table}

In a scoping study it is standard  practice to assume that the parameters are uniformly distributed allowing for an initial exploration of the sensitivities of the model \cite{CSIRO2017}. This is consistent with earlier studies by Fajraoui et.al. \cite{Fajraoui2011} for the study of porous media and  by Babaei, Alkhatib and Pan \cite{Babaei2015a} in the study of subsurface flows.

In CMG simulations, the mass conservation equations are based on the following conditions: there are only two phases (water and gas) in the coal seam; the temperature is constant; coal cleats are isotropically distributed; the gas desorption obeys the Langmuir equation; and the gas is not dissolved in water \cite{Zhou2012}. In this study, the coal shrinkage \cite{Palmer1998,Shi2005} after gas desorption is not considered.

\section{Construction of PCE surrogates}\label{sec:pce}

\subsection{The method of Polynomial Chaos Expansion}
\label{subsec:pce}
Starting with the commercial solver detailed in the previous section, we denote the outputs of the solver (cumulative gas production and peak gas extraction) by $\mathbf{Y}$, the $N$ rescaled uncertain input variables (fracture porosity, fracture permeability, reciprocal of the Langmuir pressure and Langmuir volume as listed in Table \ref{tab:inputs}) by $\boldsymbol{\xi}= (\xi_1,\dots, \xi_N)$ (so $N=4$ in this case),  and the action of the model by $\cM$. That is,
\begin{equation}
\mathbf{Y} = \cM(\boldsymbol{\xi}).
\end{equation}
A PCE  seeks to approximate the action of this model, $\cM(\boldsymbol{\xi})$, using a mathematical function of an explicitly known form $\cP(\boldsymbol{\xi})$, that is
\begin{equation}
{\cal M}(\boldsymbol{\xi}) \approx {\cal P}(\boldsymbol{\xi}).
\end{equation}
The function $\mathcal{P}(\boldsymbol{\xi})$ essentially becomes a \textit{surrogate model}, which can be used in place of the original model, $\mathcal{M}(\boldsymbol{\xi})$, provided the approximation has sufficient accuracy, as can be tested by comparing the outputs of $\mathcal{P}$ to those of $\mathcal{M}$ at `test points' sampled across the parameter space.
This is appealing because the surrogate model is a mathematical formula that can be evaluated very cheaply in comparison to executing the true model $\cM$, which requires running a computationally expensive commercial solver.

The remainder of this section is devoted to a brief general explanation of the construction of these surrogate models. Readers interested only in the specific application of the method are directed to Section \ref{sec:pce_here}.

For random variables, with finite variance, defined on bounded sets of $\mathbb{R}^n$, PCE makes use of the mathematical principle that any  function can be approximated in mean square by a sum over a basis of orthogonal polynomials. In the case where the random variable is smooth, then the approximation converges uniformly to $\cM(\boldsymbol{\xi})$ as the maximum order of the approximating orthogonal polynomials increases. The function is approximated as
\begin{equation}
\label{pce}
\cM(\boldsymbol{\xi}) \approx \cP(\boldsymbol{\xi}) = \sum_{\mathbf{i} \in \cN} \left( \mathbf{Y}_{\mathbf{i}} \prod_{j=1}^{N} P_{i_j}(\xi_j)\right),
\end{equation}
where $P_{i_j}$ is the polynomial from the orthogonal basis of degree $i_j$, in the $j$-th variable. The coefficients $\mathbf{Y}_\mathbf{i}$ need to be determined in order to construct the PCE surrogate; they are sometimes referred to as the `node strengths'.
The combination of orders for these basis polynomials is denoted by the `multi-index' $\mathbf{i}=(i_1,\dots, i_N)$, with sums taken over all possible combinations which are members of the neighbourhood $\cN$.

The neighbourhood $\cN$ can take a variety of forms,
$
\cN = \{\mathbf{i}\in\mathbb{N}^N \mid 0 \leq |\mathbf{i}|_1 \leq p\},
$
where $|\mathbf{i}|_1=\sum_{j=1}^N i_j$, is termed a \textit{total order}, while
$
\mathcal{N}  = \{\mathbf{i}\in\mathbb{N}^N \mid 0\leq i_j\leq p, j = 1,\ldots, N\}
$
is termed a \textit{tensor product}.
Which neighbourhood is used depends on how the PCE is constructed, with a full explanation given in Section \ref{sec:quadrature}.

In both of the above choices for the neighbourhood $\mathcal{N}$ we refer to $p$ as the \textit{order} of $\mathcal{N}$ and as the \textit{order} of the corresponding PCE. As $p$ is increased, more polynomials of higher order are used in the construction and thus a larger range of `shapes' can be more accurately resolved using combinations of these polynomials. Of course, as more terms are added more node strengths (coefficients) must be calculated. For convergence of the PCE, the node strengths of the higher degree polynomial terms should be small compared to earlier terms.

The statistical distributions of the unknown input variables are incorporated into the PCE by prescribing a polynomial basis, where the polynomials are orthogonal with respect to these distributions.
Some of the common bases and corresponding statistical distributions are given in Table \ref{tab:polynomial_choices}. However, orthogonal polynomial bases can be constructed for arbitrarily distributed input variables, for example via Gram-Schmidt orthogonalization, see \cite[Chapter 3]{StoerBulirsch}.

\begin{table}[ht]
\begin{center}
\resizebox{\linewidth}{!}{\begin{tabular}{llll}
\hline
\bf{Distribution} & \bf{Orthogonal Basis,} ${\cal P}(\boldsymbol{\xi})$ & \bf{Range} & \bf{Weight,} $\omega(\xi)$ \\
\hline
Uniform & Legendre Polynomials & $[-1, 1]$ & $1/2$ \\
Normal & Hermite Polynomials & $(-\infty, \infty)$ & $e^{-\xi^2/2}$\\
Exponential & Laguerre Polynomials & $[0, \infty)$ & $e^{-\xi}$ \\
Gamma & Generalised Laguerre Polynomials & $[0, \infty)$ & $\xi^{\alpha} e^{-\xi}$ \\
\hline
\end{tabular}
}
\caption{Statistical distributions aligned with associated orthogonal polynomial bases and weight functions $\omega(\xi)$.}
\label{tab:polynomial_choices}
\end{center}
\end{table}

Once the node strengths $\mathbf{Y}_{\mathbf{i}}$ are determined, Equation \eqref{pce} completely defines the PCE used to approximate the original model. Using the properties of orthogonal functions the values of these node strengths can be evaluated using the formula
\begin{equation}
\label{node_strengths}
\mathbf{Y}_{\mathbf{i}} = \frac{1}{\prod_{j=1}^{N} \langle P_{i_j}, P_{i_j} \rangle} \idotsint {\cal P}(\boldsymbol{\xi}) \prod_{j=1}^{N} P_{i_j}(\xi_j) \, d\boldsymbol{\xi},
\end{equation}
where $\prod_{j=1}^{N} \langle P_{i_j}, P_{i_j} \rangle$ is the product of the expectations of the given polynomials. We use numerical integration techniques to approximate this integral. These techniques require the values of the integrand at certain points, termed \textit{quadrature} points. Although the function ${\cal P}$ is not yet known, its value at the quadrature points can be approximated by the evaluation of the original model ${\cal M}$ at these points. This approach is known as \textit{non-intrusive} polynomial chaos, because the PCE is constructed without modification to the underlying model, albeit at the cost of running the model ${\cal M}$ at enough quadrature points to accurately approximate the integral. The quadrature rules utilised in this work are discussed in Subsection \ref{sec:quadrature}.

For the sake of illumination, we present the form of the PCE for the CMG model under discussion. In this case, there are four uncertain variables that are each assumed to take a uniform distribution. Following Table \ref{tab:polynomial_choices}, the orthogonal polynomials for each random variable are the Legendre polynomials denoted $L_i$, generated using the recurrence relation
\begin{eqnarray*}
\nonumber
L_0(\xi) &=& 1 \\
\label{legendre_polynomials}
L_1(\xi) &=& \xi \\
\nonumber
(n + 1) L_{n+1}(\xi) &=& (2n+1) \xi L_{n}(\xi) - n L_{n-1}(\xi).
\end{eqnarray*}
The random variables $\xi_1$, $\xi_2$, $\xi_3$, and $\xi_4$ are those listed in Table \ref{tab:inputs} rescaled to vary over $[-1,1]$, so that the orthogonal property of the Legendre polynomials,
\begin{equation}
\label{legendre_orthogonality}
\langle L_i, L_j \rangle = \int_{-1}^{1} L_{i}(\xi) L_{j}(\xi)\omega(\xi) \, d\xi = \left\{
\begin{array}{ll} 0,&  \mbox{for } i \neq j, \\
\dfrac{1}{2i + 1},&\mbox{for } i = j, \end{array}\right.
\end{equation}
can be utilised (recall that for uniform distributions, and hence Legendre polynomials, $\omega(\xi)=1/2$).

Restating Equation \eqref{pce} for a PCE based on uniform distributions:
\begin{equation}
\label{pce_here}
\mathbf{Y} \approx {\cal P}(\boldsymbol{\xi}) = \sum_{\mathbf{i}\in\mathcal{N}} \left(\mathbf{Y}_{\mathbf{i}} \prod_{j=1}^N L_{i_j}(\xi_j)\right).
\end{equation}
So, for instance, if  $\mathcal{N}$ is a total order neighbourhood the summation can be written as
\begin{equation}
\nonumber
{\cal P}(\xi_1,\xi_2,\xi_3,\xi_4) = \sum_{i_1 = 0}^{p} \, \, \sum_{i_2 = 0}^{p-i_1} \, \, \sum_{i_3 = 0}^{p-i_1-i_2} \, \, \sum_{i_4 = 0}^{p-i_1-i_2-i_3} \mathbf{Y}_{\mathbf{i}} L_{i_1}(\xi_1) L_{i_2}(\xi_2) L_{i_2} (\xi_3) L_{i_3} (\xi_4).
\end{equation}
The limits on the sums are used to ensure that the sum of $i_1$, $i_2$, $i_3$ and $i_4$ does not exceed $p$.

Given independent uniform distributions on $[-1,1]$ for the input variables, $\omega(\xi_j) = 1/2$, and Equation \eqref{node_strengths} yields
\begin{eqnarray}
\label{node_strengths_here}
\mathbf{Y}_\mathbf{i} &=& \prod_{j=1}^N(2i_j + 1) \times \displaystyle{\int_{[-1,1]^N} }\cP(\boldsymbol{\xi})\left(\prod_{j=1}^N L_{i_j}(\xi_j)\,\omega(\xi_j)\right)\,d\boldsymbol{\xi} \nonumber\\
&=& \prod_{j=1}^N\frac{2i_j + 1}{2}\times \displaystyle{\int_{[-1,1]^N}}\cP(\boldsymbol{\xi})\left(\prod_{j=1}^N L_{i_j}(\xi_j)\right)\,d\boldsymbol{\xi}.
\end{eqnarray}

Recalling that we approximate the $\cP(\boldsymbol{\xi})$ by $\cM(\boldsymbol{\xi})$,
\begin{equation}
\label{node_strengths2_here}
\mathbf{Y}_\mathbf{i} \approx \prod_{j=1}^N\frac{2i_j + 1}{2} \times  \displaystyle{\int_{[-1,1]^N}}\cM(\boldsymbol{\xi})\left(\prod_{j=1}^N L_{i_j}(\xi_j)\right)\,d\boldsymbol{\xi}.
\end{equation}
Hence, for our four-dimensional example 
\begin{multline}
\mathbf{Y}_\mathbf{i} \approx \frac{(2{i_1}+1)(2{i_2}+1)(2{i_3}+1)(2{i_4}+1)}{16} \times \\
\int_{[-1,1]^4} \cM(\boldsymbol{\xi}) L_{i_1}(\xi_1) L_{i_2}(\xi_2) L_{i_3} (\xi_3) L_{i_4} (\xi_4) \, d\boldsymbol{\xi},
\nonumber
\end{multline}
and we perform this integration using the quadrature techniques discussed in the following section.

\subsection{Quadrature for numerical integration}
\label{sec:quadrature}
Our goal is to accurately approximate, for each $\mathbf{i}\in\cN$,
\[
\left(\prod_{j=1}^N\frac{2i_j+1}{2}\right)\times \int_{[-1,1]^N}\cP(\boldsymbol{\xi})\left(\prod_{j=1}^N L_{i_j}(\xi_j)\right)\, d\boldsymbol{\xi}.
\]
By construction, for each summand in $\cP$, the combination of orders is an element of $\cN$. The combination of orders for each term in $\prod_{j=1}^N L_{i_j}(\xi_j)$ is also an element of $\cN$. So, if $\cN$ is a tensor product (respectively a total order) of order $p$, this implies that $\cP(\boldsymbol{\xi})\left(\prod_{j=1}^N L_{i_j}(\xi_j)\right)$ is a polynomial in which the combination of orders for any term is an element of a tensor product (respectively a total order) of order $2p$.

In each of our two cases (tensor product or total order neighbourhoods) we choose quadrature rules that are known to integrate the relevant polynomial exactly. Note however that because $\cM(\boldsymbol{\xi})$ is used  to approximate $\cP(\boldsymbol{\xi})$ in our calculation these numerical integrations may not be exact; this contributes to the overall performance which we quantify in Subsection \ref{subsec:error}.

Numerical quadrature rules approximate a multi-dimensional integral of a function, $f(\boldsymbol{\xi})$ say, as a weighted sum over values of the function at carefully selected quadrature points $\{\boldsymbol{\xi}_\mathbf{i}\mid \mathbf{i}\in \cI\}$, i.e.
\[
\int_\Omega f(\boldsymbol{\xi})\,d\boldsymbol{\xi} \approx \sum_{\mathbf{i}\in\cI} f(\boldsymbol{\xi}_\mathbf{i})w_\mathbf{i},
\]
where $w_\mathbf{i}$ is the \textit{weight} of the quadrature point $\boldsymbol{\xi}_\mathbf{i}$ and is given by the quadrature rule.

\subsubsection{Full grid quadrature}
Quadrature rules for multi-dimensional integrals can be constructed by taking tensor products of one-dimensional rules; such an approach is called \textit{full grid quadrature}. In the case where the one-dimensional rule is Gauss-Legendre quadrature $\cI =\{\mathbf{i}\in\mathbb{N} \mid 1\leq i_j\leq p+1, j = 1,\ldots,N\}$, $w_\mathbf{i} = \prod_{j=1}^N w_{i_j}$ is the product of one-dimensional Gauss-Legendre quadrature\break weights and each $\boldsymbol{\xi}_\mathbf{i}$ is the point on the grid of quadrature points with co-ordinates indexed by $\mathbf{i}$. Tables of the coordinates in one-dimension that make up this grid, and the weight coefficients associated with them, are readily available \cite{zwillinger} or can be calculated relatively simply.

It is well known, see for example \cite[Chapter 3]{StoerBulirsch}, that in one dimension Gauss-Legendre quadrature with $p+1$ quadrature points integrates any polynomial of degree less than or equal to $2p+1$ exactly. It follows by \cite[Corollary 2]{HammerWymore} that full grid Gauss-Legendre quadrature will evaluate integrals of polynomials in which the combination of orders of each term is an element of the tensor product of order $2p+1$ exactly. 

As the components of each vector in $\cI$ can take integer values between $1$ and $p+1$ inclusively, this results in $(p+1)^N$ quadrature points being required to compute the integral. This number grows very quickly as $N$ increases and hence this approach is computationally costly when the number of dimensions (uncertain input variables) is large.

\subsubsection{Sparse grid quadrature}
The above limitation of full grid approaches can be somewhat addressed through a procedure due to Smolyak \cite{smolyak_quadrature}: \textit{sparse grid quadrature}. The resulting approximation of the integral $A(\ell,N)\approx\int_{[-1,1]^N} f(\boldsymbol{\xi})\,d\boldsymbol{\xi}$ can be expressed as 
\begin{eqnarray}
A(\ell,N) &=& \sum_{\ell+1 \leq |\mathbf{k}|_1 \leq \ell+N} (-1)^{\ell + N - |\mathbf{k}|_1}\binom{N-1}{\ell+N-|\mathbf{k}|_1}Q_\mathbf{k},
\label{sparsegridquad}
\end{eqnarray}
where $Q_\mathbf{k}=Q_1^{k_1}\otimes\ldots\otimes Q_N^{k_N}$ and $Q^{k_j}_j$ is a one-dimensional quadrature rule that exactly integrates all univariate polynomials of degree $k_j$. In the case where the one-dimensional quadrature rules have symmetric and bounded weight functions and $\ell \leq 3N$, $A(\ell,N)$ exactly computes integrals of polynomials where the combination of orders is an element of the total order neighbourhood of order $2\ell+1$ \cite[Corollary 3]{NovakRitter}; we call $\ell$ the level of the sparse grid quadrature. We use Clenshaw-Curtis quadrature rules to construct our sparse grid quadrature, which together with the choice for the number of points $n_{k_j} = 2^{k_j-1} + 1$ provides the important property that the set of quadrature points required for the one-dimensional quadrature rule of level $k$ contains the points of all levels up to $k$; i.e., they are \textit{nested} quadrature rules. This means that the quadrature rules of the different levels up to $\ell$ that are summed in Equation \eqref{sparsegridquad} contain the quadrature points of the previous level, greatly reducing the total number of model evaluations required. The sparse grid approach is non-trivial to implement, however as $N$ increases the number of quadrature points required grows far more slowly than in the full grid case.

Note that the nested property is an advantage if the degree of the polynomial to be integrated is increased, the set of quadrature points for the new order contains all the quadrature points for lower degree polynomials.  Hence the overall number of expensive evaluations of $\cM(\boldsymbol{\xi})$ is reduced. For example, suppose that PCEs are to be constructed for a model of dimension four initially at order 4 and then, if the order 4 PCE does not sufficiently approximate the original model, at order 5. In the case of full-grid Gauss-Legendre quadrature we would require $625+1296=1931$ quadrature points (evaluations of $\cM(\boldsymbol{\xi})$) while for sparse grid  Clenshaw-Curtis quadrature we require $1105$ quadrature points (evaluations of $\cM(\boldsymbol{\xi})$) as the $401$ points for $p = 4$ are repeated in the $1105$ points for $p = 5$ (Smolyak's procedure preserves the nested property, see \cite{NovakRitter1996}).

\section{Global sensitivity analysis}
\label{subsubsec:GSA}

As we shall shortly formalise, a PCE can be utilised to perform a global sensitivity analysis for a model at negligible computational cost.
The aim of a global sensitivity analysis of a model is to apportion the variance in the output of the model to the variances of the inputs or combinations of inputs. Often in order to undertake such an analysis many function evaluations must be made, however, through the use of the PCE surrogate these expensive evaluations can be avoided.

We now provide some background to a method due to Sobol' \cite{Sob_Rus, Sob_Eng}, also see \cite{SalChaSco}, to undertake global sensitivity analysis of  a model $\cL$ of dimension (number of variables) $N$ where each variable is drawn from a closed interval.
As, we may rescale each of the model inputs, we may assume that the domain of the model $\cL$ is $[-1,1]^N$. We assume that the model inputs $\boldsymbol{\xi}$ are independent and denote the joint pdf as $\omega(\boldsymbol{\xi})$. We will also use the notation $[N]=\{1, \ldots, N\}$.

Suppose that $U=\{u_1, u_2,\ldots, u_k\}$ is a subset of $[N]$ of cardinality $k$ and assume, without loss of generality, that $u_1\leq u_2\leq ...\leq u_k$. Then for  $\boldsymbol{\xi}=(\xi_1,\xi_2,\ldots,\xi_N)\in [-1,1]^N$ we denote the vector $(\xi_{u_1}, \xi_{u_2}, \ldots, \xi_{u_k})$ by $\boldsymbol{\xi}_U$ and the joint pdf of these variables %$\xi_{u_1}, \xi_{u_2}, \ldots, \xi_{u_k}$
 as $\omega(\boldsymbol{\xi}_U)$.

Sobol's method is based on expressing $\cL$ as a decomposition into summands of increasing dimensionality, that is
\begin{equation}
\cL(\boldsymbol{\xi}) = m_0 +  \sum_{\alpha=1}^N\left(\sum_{U\subseteq [N], |U|=\alpha} m_U(\boldsymbol{\xi}_U)\right),
\label{eqn:decomp_into_summands}
\end{equation}
where $m_U(\boldsymbol{\xi}_U)$ is a function satisfying  $\int_{-1}^1 m_U(\boldsymbol{\xi}_U)\,\omega(\boldsymbol{\xi}_U)\, d{\xi_u} = 0$ for all $u\in U$ and  $m_0$ is a constant.
It follows that the summands are orthogonal; i.e. for $U,V\subseteq [N]$ and $U\neq V$ we have
$
\int_{[-1,1]^N} m_U(\boldsymbol{\xi}_U)\,m_V(\boldsymbol{\xi}_V)\,\omega(\boldsymbol{\xi})\,d\boldsymbol{\xi} = 0.
$
Sobol' \cite{Sob_Eng}, also see \cite{SalChaSco},  showed that for $U\subseteq [N]$ the summands can be calculated as follows.
\begin{eqnarray*}
m_0 & = & \int_{[-1,1]^N} \cL(\boldsymbol{\xi})\omega(\boldsymbol{\xi})\,d\boldsymbol{\xi};\text{ and}\\
m_U(\boldsymbol{\xi}_U) & = & -m_0  - \sum_{V\subset U}m_V(\boldsymbol{\xi}_V)  +\int_{[-1,1]^{N-|U|}} \cL(\boldsymbol{\xi})\omega(\boldsymbol{\xi}_{[N]\setminus U})\,d\boldsymbol{\xi}_{[N]\setminus U},
\end{eqnarray*}
where $m_\emptyset(\boldsymbol{\xi}_\emptyset)$ is interpreted as $0$. Note that $m_0$ is the mean of $\cL$.

Given this decomposition the variance of $\cL$ can be expressed as:
\begin{equation}
\varepsilon(\cL) = \int_{\Omega^N} \cL^2(\boldsymbol{\xi})\omega(\boldsymbol{\xi})\,d\boldsymbol{\xi} - (m_0)^2 =\sum_{i=1}^N\left(\sum_{U\subseteq[N],|U|=i}D_U\right),
\label{eqn:decomp_var}
\end{equation}
where
\[
D_U = \int_{\Omega^{|U|}} m^2_U(\boldsymbol{\xi}_U)\,d\boldsymbol{\xi}_U.
\]
That is, $D_U$ is the portion of the variance due to the input variables $\boldsymbol{\xi}_U$.

The \textit{Sobol'} or \textit{variance-based sensitivity indices}  are obtained by normalising these in terms of the total variance, that is
\[
S_{\boldsymbol{\xi}_U} = \frac{D_U}{\varepsilon(\cL)}.
\]
Hence these indices express the fractional contribution to the variance of $\cM$ from each input or combination of inputs. If $|U|=1$, then the index is referred to as a \textit{main effect Sobol' index}.
The \textit{total Sobol' index} of a set of inputs $\boldsymbol{\xi}_U$ is the sum of all Sobol' indices whose indexing sets involve all the elements of $\boldsymbol{\xi}_U$, see \cite{SalChaSco}, we denote this index by $T_{\boldsymbol{\xi}_U}$, that is
\[
T_{\boldsymbol{\xi}_U} = \sum_{V\supseteq U} S_{\boldsymbol{\xi}_V}.
\]

\begin{example}
Let $\xi_1$ and $\xi_2$ be independent uniform random variables each drawn from $[-1,1]$ and consider the models (functions)
\[
\cL_1(\xi_1,\xi_2) = \xi_1^2 + \xi_2^2\quad\text{ and }\quad \cL_2(\xi_1,\xi_2)  = \xi_1^3 + \xi_2.
\]
Note that the pdfs for  $\xi_1$ and $\xi_2$ are $\omega(\xi_1) = \omega(\xi_2) = 1/2$ and their joint pdf is $\omega(\boldsymbol{\xi}) = 1/4$.

Following the above, the decompositions of $\cL(\xi_1,\xi_2)$ and $\cL(\xi_1,\xi_2)$ in the form shown in Equation (\ref{eqn:decomp_into_summands}) are:
\[
\cL_1(\xi_1,\xi_2)  = \frac{2}{3} + \left(-\frac{1}{3} + \xi_1^2\right) + \left(-\frac{1}{3} +\xi_2^2\right) + 0
\]
and
\[
\cL_2(\xi_1,\xi_2)  = 0+ \xi_1^3 + \xi_2 + 0.
\]
We now calculate their variances and the terms of their decompositions as given in Equation (\ref{eqn:decomp_var}):
\begin{eqnarray*}
\varepsilon(\cL_1) &=& \int_{-1}^1\int_{-1}^1\left(\xi_1^2+\xi_2^2\right)^2\, \frac{1}{4}\,d\boldsymbol{\xi} - \left(\frac{2}{3}\right)^2 = \frac{8}{45},\\
D^{(1)}_{\{1\}} &=& \int_{-1}^1 \left( -\frac{1}{3} + \xi_1^2\right)^2\,\frac{1}{2}\, d\xi_1 = \frac{4}{45},\\
D^{(1)}_{\{2\}} &=& \int_{-1}^1 \left( -\frac{1}{3} + \xi_2^2\right)^2\,\frac{1}{2}\, d\xi_2 = \frac{4}{45},\\
\varepsilon(\cL_2) &=& \int_{-1}^1\int_{-1}^1\left(\xi_1^3+\xi_2\right)^2\, \frac{1}{4}\,d\boldsymbol{\xi} - 0 = \frac{10}{21},\\
D^{(2)}_{\{1\}} &=& \int_{-1}^1 \xi_1^6\,\frac{1}{2}\,d\xi_1 = \frac{1}{7},\text{ and} \\
D^{(2)}_{\{2\}} &=& \int_{-1}^1 \xi_2^2\,\frac{1}{2}\, d\xi_2 = \frac{1}{3},
\end{eqnarray*}
where the bracketed superscripts indicate the respective models.
Hence, the Sobol' indices for $\cL_1$ are $S^{(1)}_{(\xi_1)} = 1/2$ and $S^{(1)}_{(\xi_2)} = 1/2$;
while the Sobol' indices for $\cL_2$ are $S^{(2)}_{(\xi_1)}= 3/10$ and $S^{(2)}_{(\xi_2)} = 7/10$.
Note that $\xi_1$ and $\xi_2$ contribute equally to the variance of $\cL_1$, while the variance of $\cL_2$ is dominated by the contribution due to $\xi_2$.
\end{example}

In general calculating the integrals to evaluate the above indices may be computationally expensive. However, given a function/model $\cP$ in the form of a PCE expansion, i.e.
\[\cP(\boldsymbol{\xi}) = \sum_{\mathbf{i} \in \cN} \left( \mathbf{Y}_{\mathbf{i}} \prod_{j=1}^{N} P_{i_j}(\xi_j)\right),\]
Sudret \cite{Sudret2008} showed that the Sobol' index for a set of variables $\boldsymbol{\xi}_U$ can be analytically calculated using the coefficients of the expansion as
\[S_{\boldsymbol{\xi}_U} = \frac{1}{\varepsilon(\cP)} \sum_{\mathbf{i}\in\mathcal{I}_U} \mathbf{Y}_\mathbf{i}^2\langle L_\mathbf{i},L_\mathbf{i}\rangle \]
where $\mathcal{I}_\alpha = \{ (i_1,i_2,\ldots,i_N) : i_j>0 \text{ if and only if } j\in\alpha, \text{ for all }1\leq j\leq N \}$. Note that, as $\langle L_\mathbf{0},L_\mathbf{0}\rangle=1$, $\varepsilon(\cP)$ can be easily calculated by
\[\varepsilon(\cP) = \sum_{\mathbf{i}\in\cN}\mathbf{Y}_\mathbf{i}^2\langle L_\mathbf{i},L_\mathbf{i}\rangle - \mathbf{Y}_\mathbf{0}^2.\]

In the case where the variables are independent and drawn uniformly from $[-1,1]$ the polynomials $L_\mathbf{i}$ are products of univariate Legendre polynomials and hence
\[S_{\boldsymbol{\xi}_U}  = \frac{1}{\varepsilon(\cP)} \sum_{\mathbf{i}\in\mathcal{I}_\alpha} \frac{\mathbf{Y}_\mathbf{i}^2}{\prod_{j=1}^N(2i_j +1) }.\]

As discussed above we can approximate a model $\cM$ with a PCE $\cP$ with convergence in the mean square sense, and hence we can approximate the Sobol' indices of $\cM$ with the Sobol' indices of $\cP$. This allows for global sensitivity analysis of the commercial model $\cM$.

\section{Surrogates for  the commercial solver}
\label{sec:pce_here}
The commercial solver studied here has four uncertain inputs, denoted $V_j$, with domains as listed in Table \ref{tab:inputs}. For this model the inputs $V_1, \ldots, V_4$ are assumed to be independent and uniformly distributed over their specified ranges and after re-scaling correspond to,
\begin{equation}
\label{rescaling}
\xi_j  = \frac{2}{V_j^{\text{max}}-V_j^\text{min}}\left(V_j - \left(V_j^{\text{min}} + \frac{V_j^{\text{max}} - V_j^{\text{min}}}{2} \right) \right) = \frac{2V_j - V_j^{\text{max}} - V_j^{\text{min}}}{V_j^{\text{max}} - V_j^{\text{min}}}.
\end{equation}
So the $\xi_j$ vary uniformly over  $[-1,1]$.
The PCE surrogate is then constructed as a sum over polynomials in terms of these rescaled variables, Equation \eqref{pce_here}, with the node strengths defined using Equation \eqref{node_strengths2_here}. The four-dimensional integral in this equation is approximated using either full grid Gauss-Legendre quadrature or sparse grid Clenshaw-Curtis quadrature depending on the choice of neighbourhood $\cN$.

To execute the chosen quadrature rule, the value of the integrand in Equation \eqref{node_strengths2_here} must be evaluated at each of the quadrature points and the integral must be re-calculated for each node strength. However, as each integrand uses the same quadrature points and the computational cost is dominated by the commercial model's evaluations, this does not significantly contribute to the overall computational cost.
Thus the cost of the PCE surrogate model is calculated as the number of commercial model evaluations and is equal to the number of quadrature points used in the construction. In commercial situations where a very large number of time consuming model simulations may be required (such as uncertainty quantification or parameter estimation), the PCE surrogate thus represents a significant saving.

\subsection{Error estimation and computational cost}
\label{subsec:error}

First, the predictions of the generated PCE surrogates were tested against the outputs of the original model, that is, the results of the commercial solver. To determine error across the entire four-dimensional parameter space, trial points were chosen using Latin hypercube sampling, a method that provides superior spread of sample points over a multidimensional parameter space \cite{McKay1979}. The method works by dividing the domain of each uncertain variable into $n$ subdivisions of equal width, and then repeatedly choosing samples of $n$ points such that exactly one of the $n$  points falls within each division along each dimension.
Here 300 samples with $n = 10$ were generated, corresponding to 3000 points and an expected coverage of 26 \% of  the parameter space \cite{Donovan2016}. We will refer to these as the Latin hypercube (LH) test points.

\begin{figure}
\begin{center}

\resizebox{\linewidth}{!}{\begin{tabular}{cc}
\includegraphics[width = 6cm]{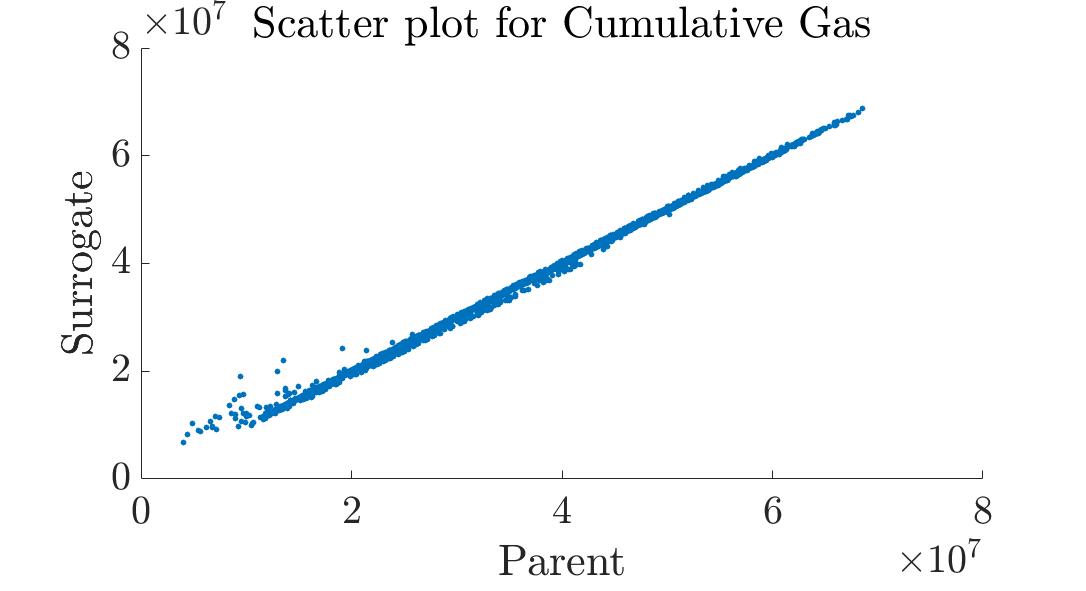} & \includegraphics[width = 6cm]{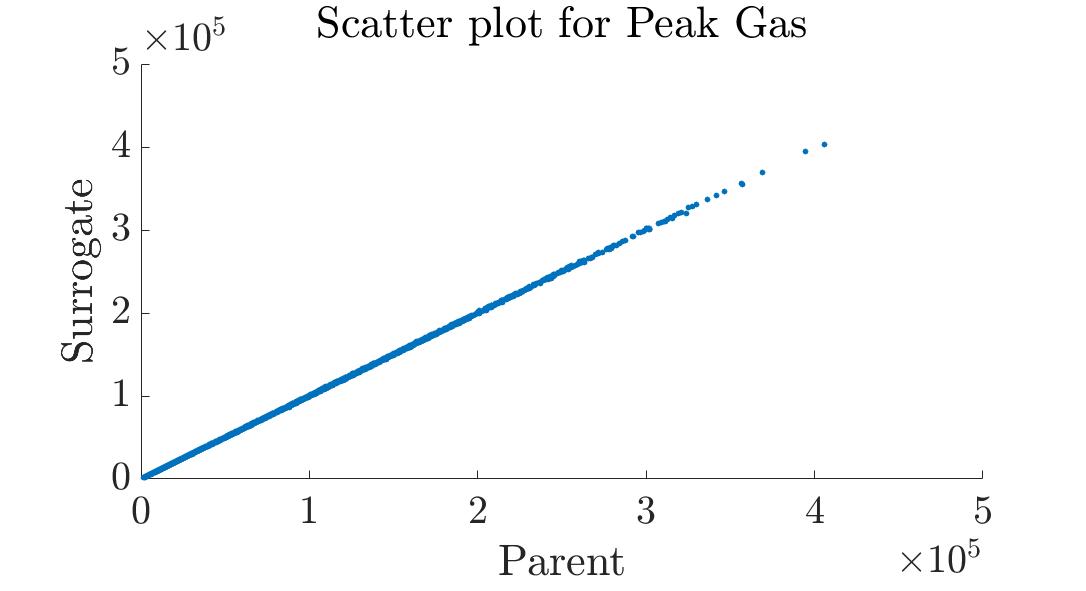}\\
\end{tabular}}

\end{center}
\caption{Scatter plots of the CMG model and the full grid Gauss-Legendre $p=6$ PCE model evaluations at the 3000 LH test points.}
\label{fig:scatter_GL_6}
\end{figure}

We use these evaluations to calculate empirical estimations for the error between the constructed surrogates and the original model.  An initial visualisation of the predictive accuracy of the surrogates can be obtained through scatter plots comparing the original model evaluations at the LH test points, see Figure \ref{fig:scatter_GL_6} for example.  The plots indicate excellent agreement between the surrogates and original models.

Because the background theory for PCE guarantees that, as $p$ increases, the expansion approaches the model in the mean square sense,  it is natural to consider error estimates based upon such a metric. Hence, in our analysis, we consider root mean square error and the dimensionless relative root mean square error.

Let $\mathcal{X}$ be a set of $n$ test points, then the root mean square error (RMSE) is calculated as
\begin{equation}
\label{rmse}
\text{RMSE} = \sqrt{\frac{1}{n} \sum_{\mathbf{x}\in\mathcal{X}} \big(\cP(\mathbf{x}) - \cM(\mathbf{x})\big)^2}.
\end{equation}
This measures how well the different PCE surrogate models approximate the model across the entire parameter space. Normalising by the model's evaluations at the points gives the dimensionless relative root mean square error (rRMSE):
\begin{equation}
\label{rrmse}
\text{rRMSE} = \sqrt{\frac{1}{n} \sum_{\mathbf{x}\in\mathcal{X}} \left(\frac{\cP(\mathbf{x}) - \cM(\mathbf{x})}{\cM(\mathbf{x})}\right)^2}.
\end{equation}

The RMSE and rRMSE for the various PCE surrogates, alongside the number of model evaluations required to construct the PCE, are shown in Tables \ref{tab:RMSE} and \ref{tab:rRMSE}.
As expected increasing the polynomial order improves the accuracy of the method in nearly all cases. %{\color{red} 
Typically as the number of quadrature points increases the accuracy of the subsequent PCEs also increases; an exception occurs for the cumulative gas models where the sparse grid on $1105$ points outperforms the full grid on $1296$ points.%}

\begin{table}[!ht]
\centering
\begin{tabular}{lrrr}%{cccc}
\hline
\bf{Quadrature \& Order} & \bf{Cumulative Gas} & \bf{Peak Gas} & \bf{Evaluations} \\
\hline
Sparse, $p=4$ & $1.41 \times 10^{6}$ & $2.08 \times 10^{3}$ & 401 \\
Sparse, $p=5$ & $0.93 \times 10^{6}$ & $1.67 \times 10^{3}$ & 1105 \\
Full, $p=5$ & $0.88 \times 10^{6}$ & $0.46 \times 10^{3}$ & 1296 \\
Full, $p=6$ & $0.58 \times 10^{6}$ & $0.38 \times 10^{3}$ & 2401 \\
\hline
\end{tabular}
\caption{Predictive performance of PCE surrogates as compared to the original model, measured by RMSE. }
\label{tab:RMSE}
\end{table}

\begin{table}[!ht]
\centering
\resizebox{\linewidth}{!}{\begin{tabular}{lrrr}
\hline
\bf{Quadrature \& Order} & \bf{Cumulative Gas} & \bf{Peak Gas} & \bf{Evaluations} \\
\hline
Sparse, $p=4$ & $9.74 \times 10^{-2}$ & $15.39 \times 10^{-2}$ & $401$ \\
Sparse, $p=5$ & $6.45 \times 10^{-2}$ & $8.66 \times 10^{-2}$ & $1105$ \\
Full, $p=5$ & $8.83 \times 10^{-2}$ & $0.67 \times 10^{-2}$ & $1296$ \\
Full, $p=6$ & $ 5.75 \times 10^{-2}$ & $0.73 \times 10^{-2}$ & $2401$ \\
\hline
\end{tabular}}
\caption{Predictive performance of PCE surrogates as compared to the original model, measured by rRMSE. }
\label{tab:rRMSE}
\end{table}

\subsection{Uncertainty quantification}

Having obtained a PCE for the model we can quickly evaluate its summary statistics and provide further quantification of uncertainties; this analysis then forms an approximation of the equivalent analysis for the original model. In the following discussion we will use as an exemplar the full grid Gauss-Legendre PCE of order 6.

Some of the summary statistics of the PCE can be evaluated analytically from the expansion (e.g. the mean and standard deviation) while others are empirical estimates obtained from evaluating the PCE at a selection of points. We base the empirical estimates on 3000 LH sample points. The evaluation of all 3000 points takes under a second (on an Intel(R) Core(TM) i7-4770 CPU at 3.40GHz), comparing favourably to evaluating the original model, in which a single evaluation takes about five minutes. The resulting summary statistics are shown in Table \ref{tab:sum_stats}. Summary statistics for the original models are shown in Table \ref{tab:OM_sum_stats}; the derivation of all of these statistics is empirical and is based on the evaluation of the original models at 3000 LH sample points (the set of LH points is different from the one chosen for the evaluations the PCEs in Table \ref{tab:sum_stats}). 
The summary statistics indicate how successfully the PCE surrogate models approximate the probability distributions of the outputs of the original models, with close agreement between mean values, standard deviations and percentiles. The discrepancies between the surrogates minimums and maximum values and those of the original models are not unexpected as the points evaluated to establish these samples are taken from two independent sets of LHS.

\begin{table}[ht]
\centering
\begin{tabular}{lrrl}
\hline
\bf{Statistic} & \bf{Cumulative Gas} & \bf{Peak Gas} & \bf{Derivation} \\
 & $\times 10^7\ m^3$ & $\times 10^5\ m^3$ &  \\
\hline
Mean & $3.469$& $1.002$& Analytic \\
Standard Deviation & $1.387$& $0.673$ & Analytic \\
Sample minimum & $0.736$ & $0.015$& Empirical \\
10th percentile & $1.657$ & $0.192$& Empirical \\
25th percentile & $2.272$ & $0.476$ & Empirical \\
50th percentile & $3.400$ & $0.916$ & Empirical \\
75th percentile & $4.601$ & $1.392$ & Empirical \\
90th percentile & $5.419$ & $1.872$ & Empirical \\
Sample maximum  & $6.836$ & $3.887$ & Empirical \\
\hline
\end{tabular}
\caption{Summary statistics from full grid Gauss-Legendre PCE of order 6. }
\label{tab:sum_stats}
\end{table}

\begin{table}[!ht]
\centering

\begin{tabular}{lrrl}
\hline
\bf{Statistic} & \bf{Cumulative Gas} & \bf{Peak Gas} & \bf{Derivation} \\
 & $\times 10^7\ m^3$ & $\times 10^5\ m^3$ &  \\
\hline
Mean & $3.462$& $0.998$& Empirical \\
Standard Deviation & $1.390$& $0.670$ & Empirical \\
Sample minimum & $0.407$ & $0.014$& Empirical \\
10th percentile & $1.644$ & $0.200$& Empirical \\
25th percentile & $2.307$ & $0.475$ & Empirical \\
50th percentile & $3.380$ & $0.910$ & Empirical \\
75th percentile & $4.590$ & $1.400$ & Empirical \\
90th percentile & $5.422$ & $1.885$ & Empirical \\
Sample maximum  & $6.854$ & $4.060$ & Empirical \\
\hline
\end{tabular}

\caption{{Summary statistics from the original models. }}
\label{tab:OM_sum_stats}
\end{table}

\begin{figure}
\begin{center}

\includegraphics[width = 9cm]{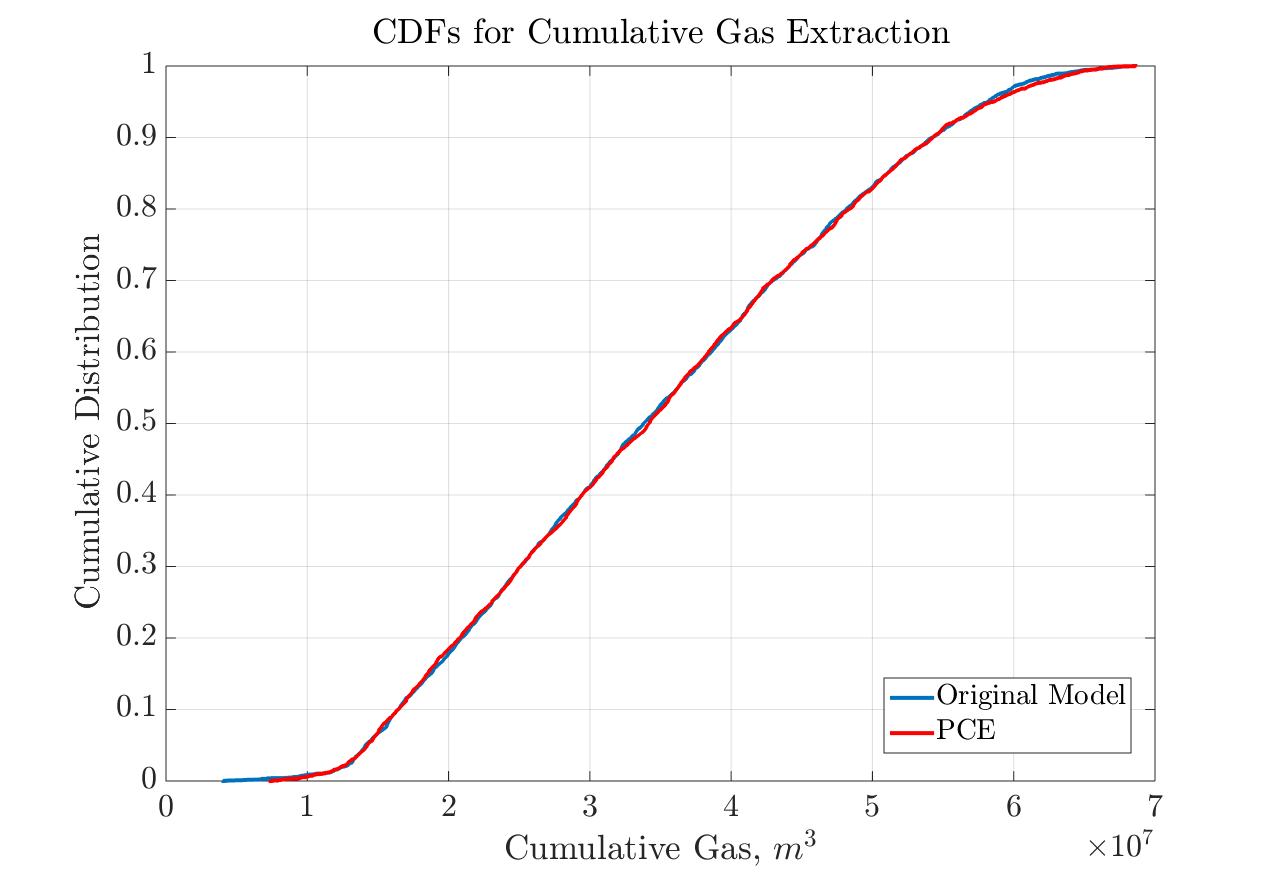} 

\vspace{10mm}

\includegraphics[width = 9cm]{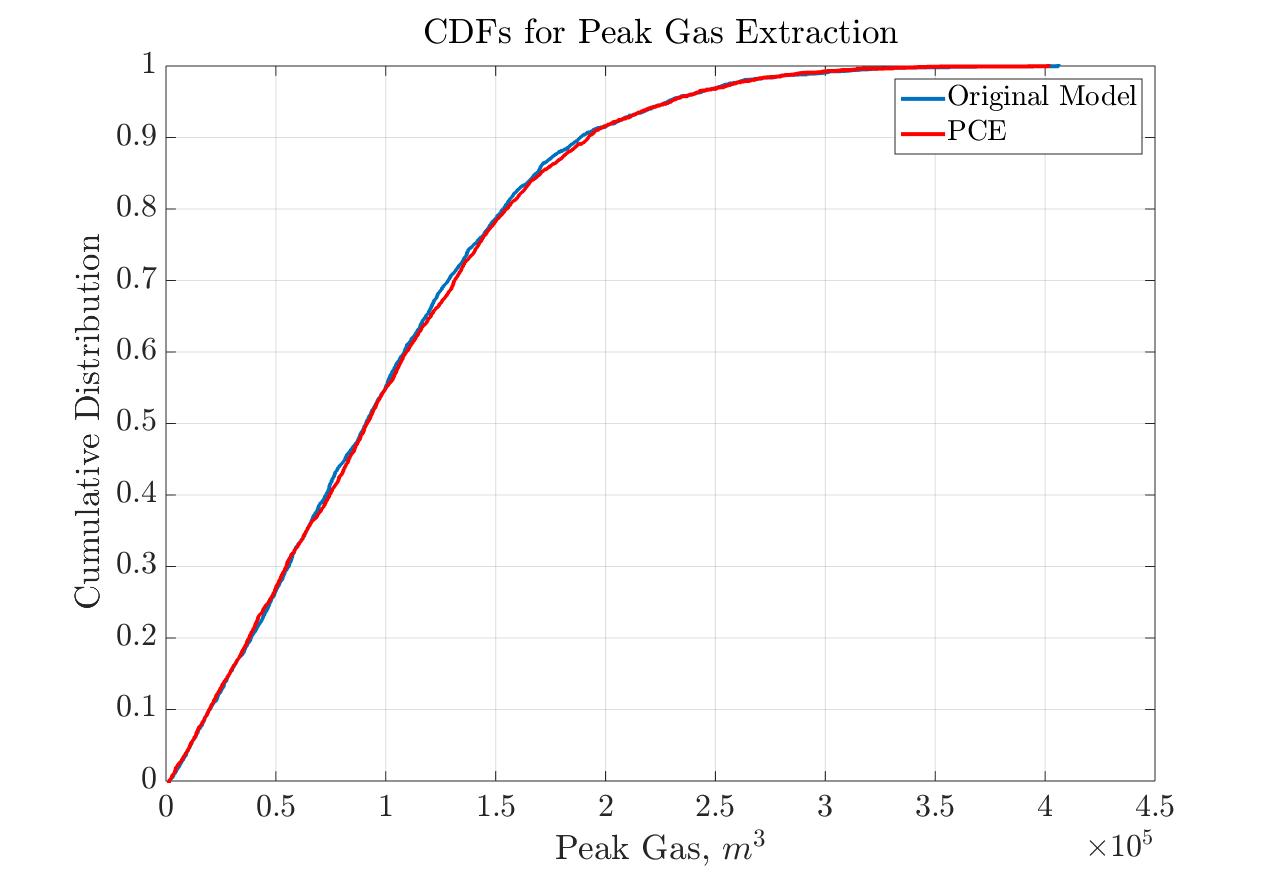}

\end{center}
\caption{Plots of empirical cumulative distribution functions of both the original models and the full grid PCEs of order six (3000 evaluations of each). }
\label{fig:cdfs}
\end{figure}

\begin{figure}
\begin{center}

\includegraphics[width = 12cm]{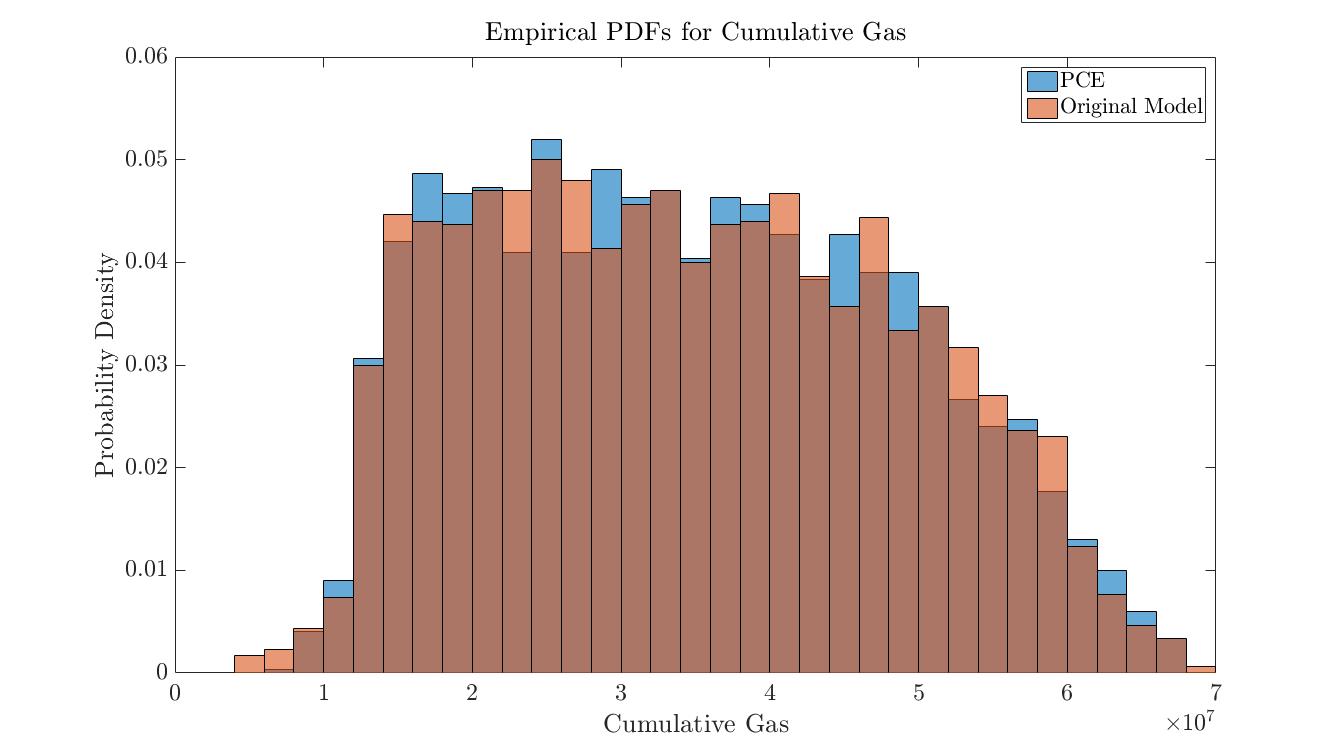} 

\vspace{10mm}

\includegraphics[width = 12cm]{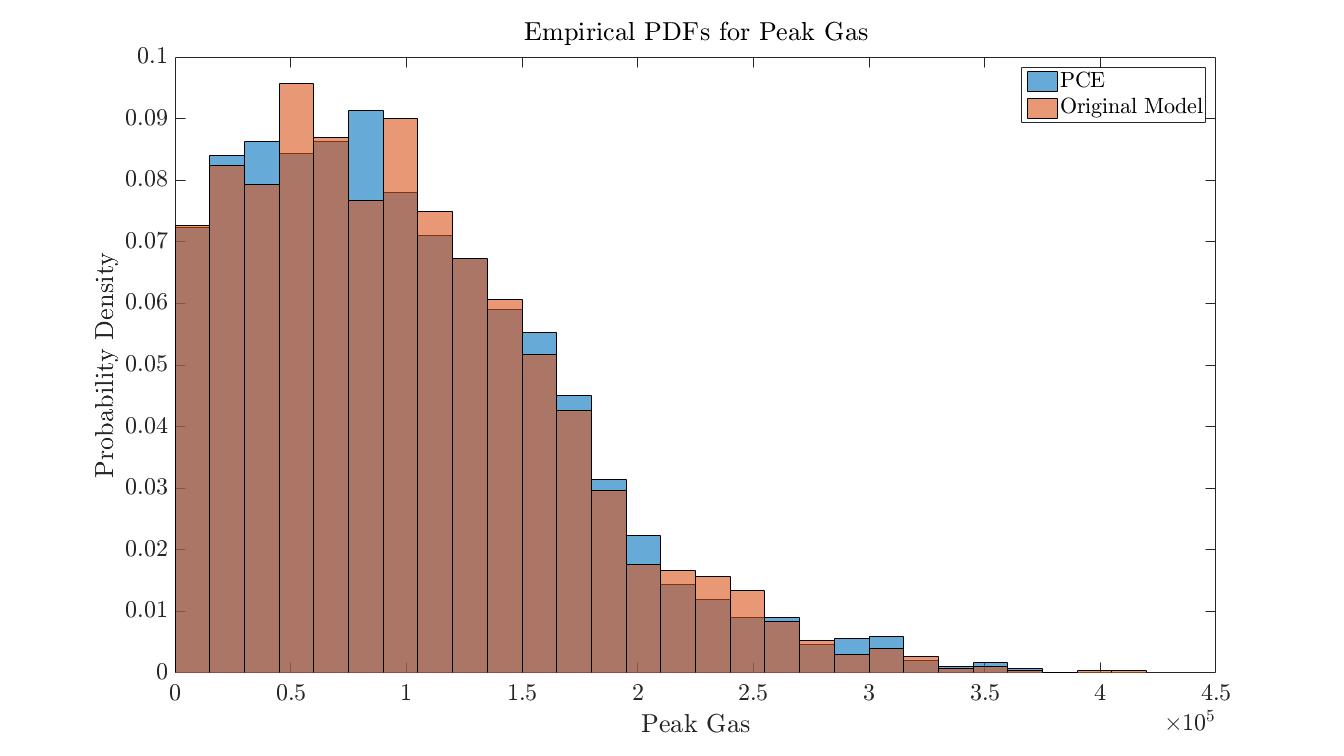}

\end{center}
\caption{
Plots of empirical probability distributions of both the original models and the full grid PCEs of order 6 (3000 evaluations of each).}
\label{fig:OM_pdfs}
\end{figure}

Empirical cumulative distribution functions for both the original models and the full grid Gauss-Legendre PCEs of order 6, all based on 3000 evaluations of the models, are shown in Figure \ref{fig:cdfs}.
 Figure \ref{fig:OM_pdfs} shows empirical probability density functions of the $3000$ evaluations of the PCEs and of $3000$ evaluations of the original models (the same evaluations used in Tables \ref{tab:sum_stats} and \ref{tab:OM_sum_stats}). These figures illustrate the success of the surrogate models in predicting the probability distributions of the outputs of the original models.

\begin{figure}
\begin{center}

\includegraphics[width = 10cm]{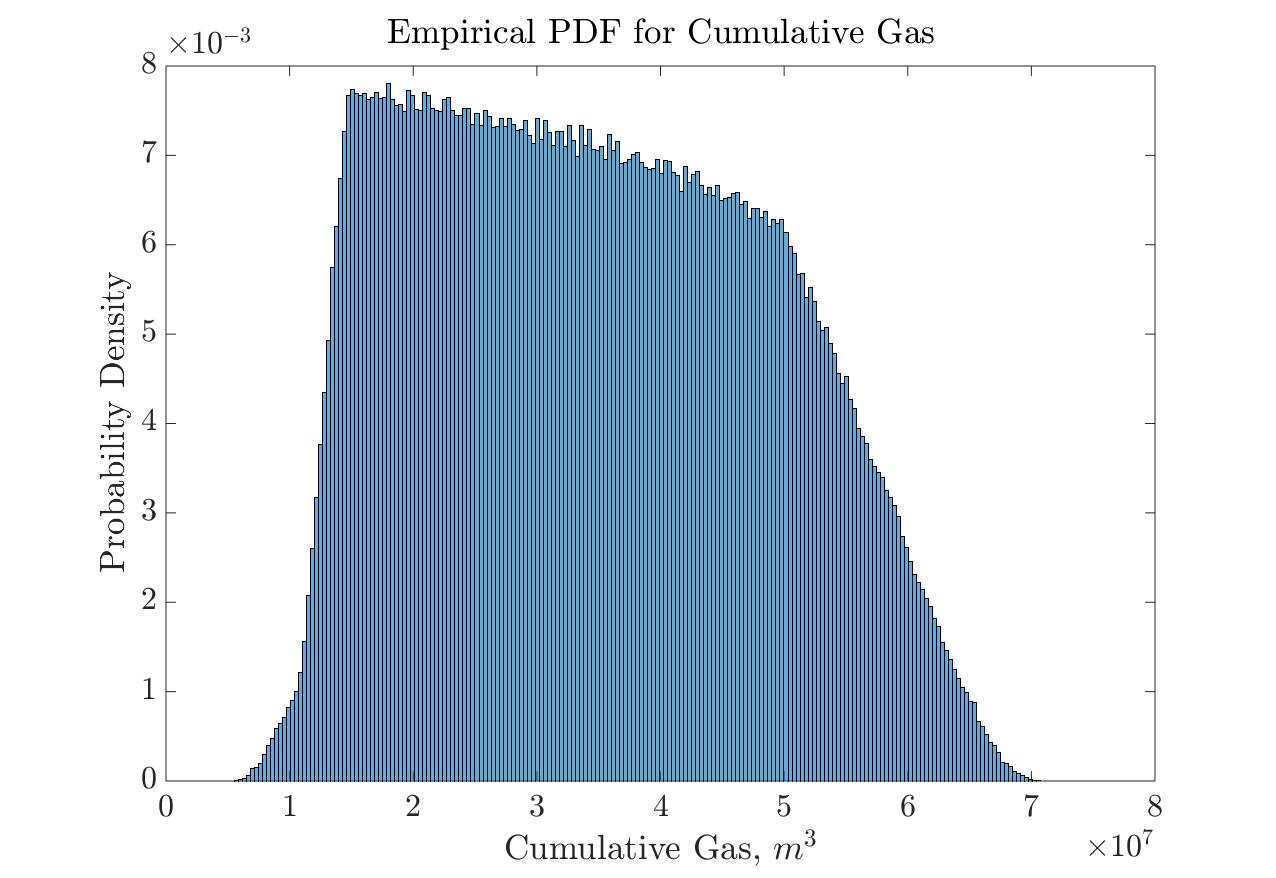} 

\vspace{10mm}

\includegraphics[width = 10cm]{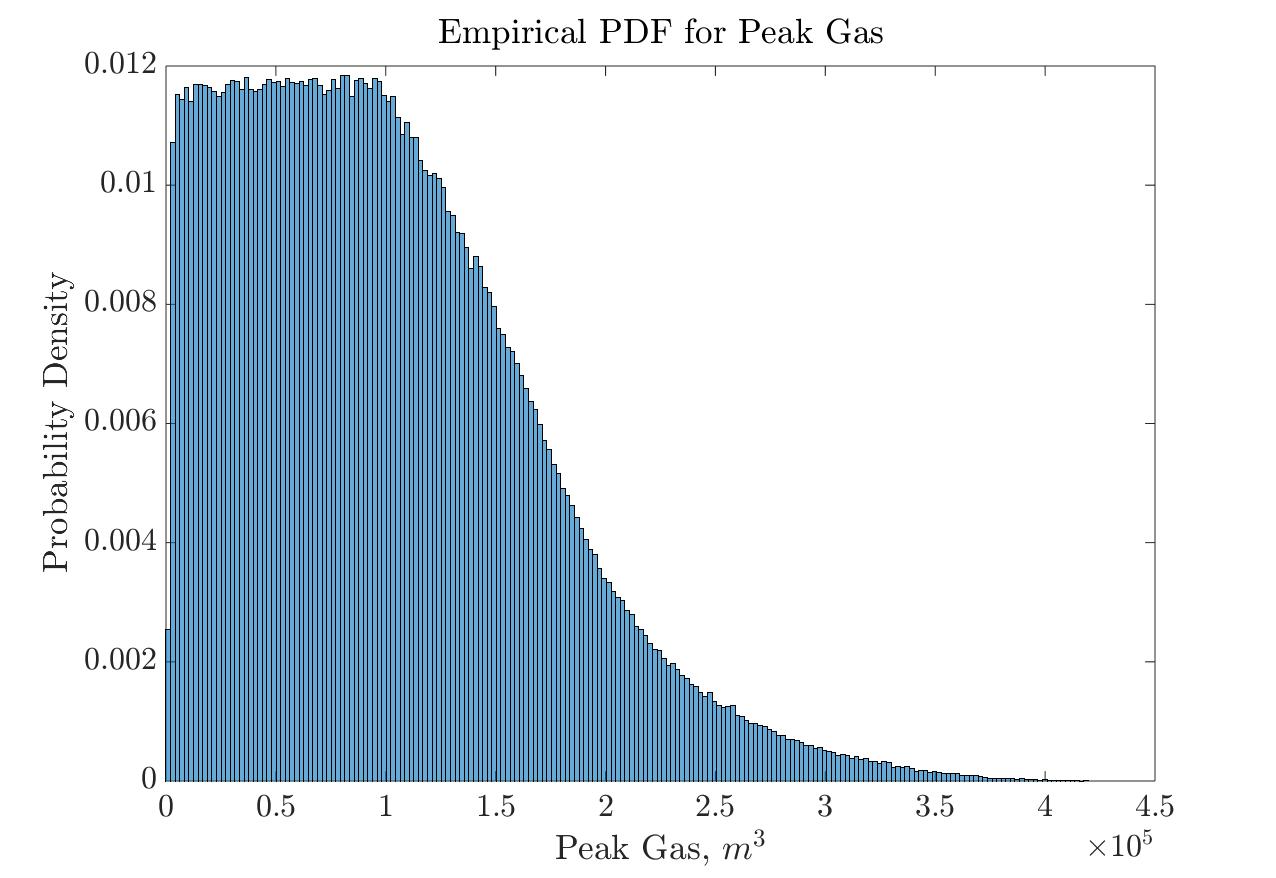}

\end{center}
\caption{Empirical probability distributions for the full grid Gauss-Legendre PCEs of order 6 (based on $1\, 000\, 000$ evaluations of the PCE).}
\label{fig:pdfs}
\end{figure}

The previous discussion focusses on the success of the PCEs to approximate the original models.
Figure \ref{fig:pdfs} shows the empirical probability distributions for the  full grid Gauss-Legendre PCEs of order 6. These plots are based on $1\, 000\, 000$ evaluations of the PCE; these evaluations took just over 203
seconds on an Intel(R) Core(TM) i7-4770 CPU at 3.40GHz. The comprehensive sampling used to construct such empirical probability density functions illustrates a benefit of the computational savings achieved through constructing the PCEs.

\subsection{Global sensitivity analysis via a surrogate}

\begin{figure}
\begin{center}

\resizebox{\linewidth}{!}{\begin{tabular}{cc}
\includegraphics[width = 6cm]{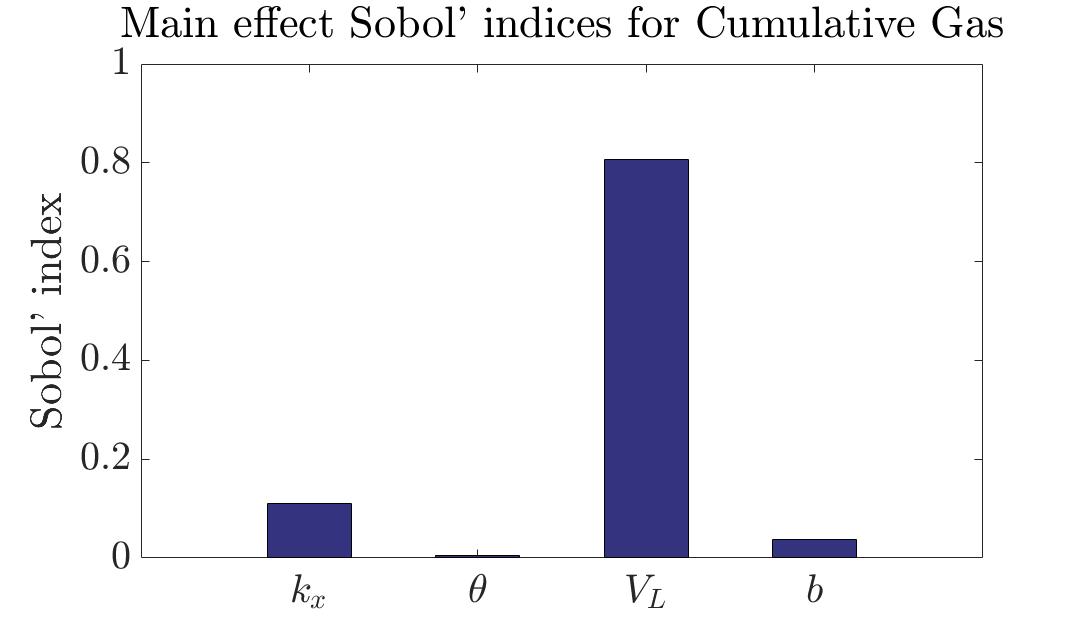} & \includegraphics[width = 6cm]{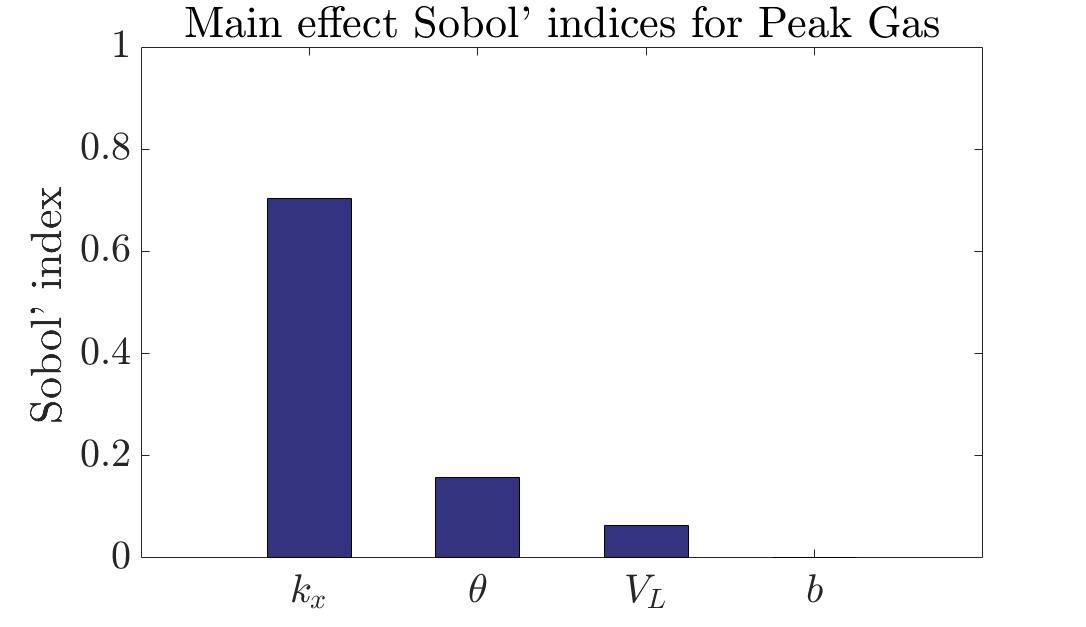}\\
&\\
\includegraphics[width = 6cm]{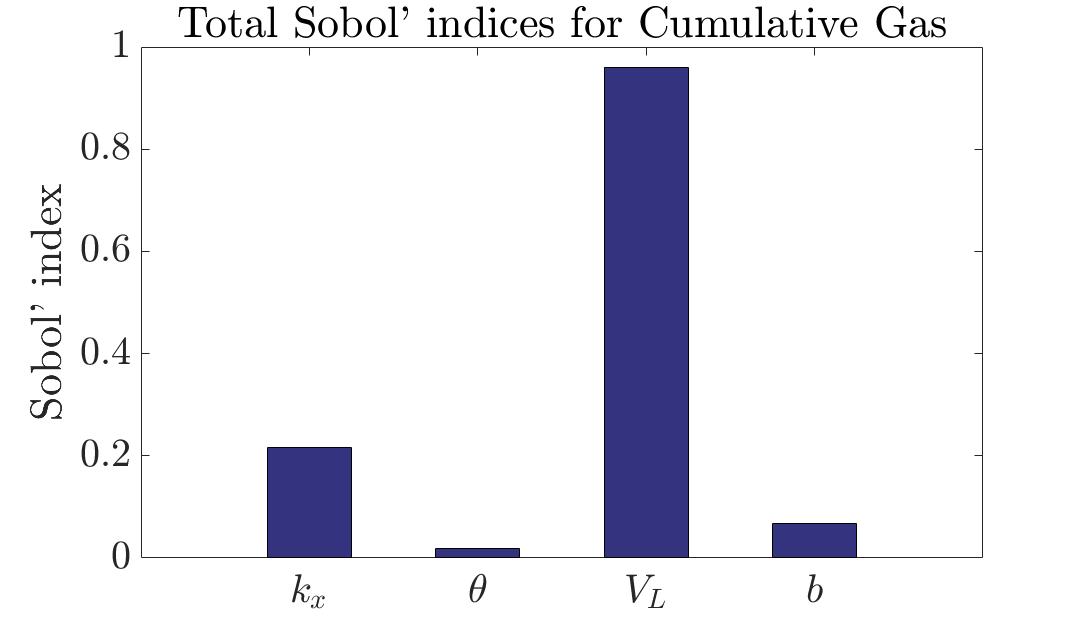} &  \includegraphics[width = 6cm]{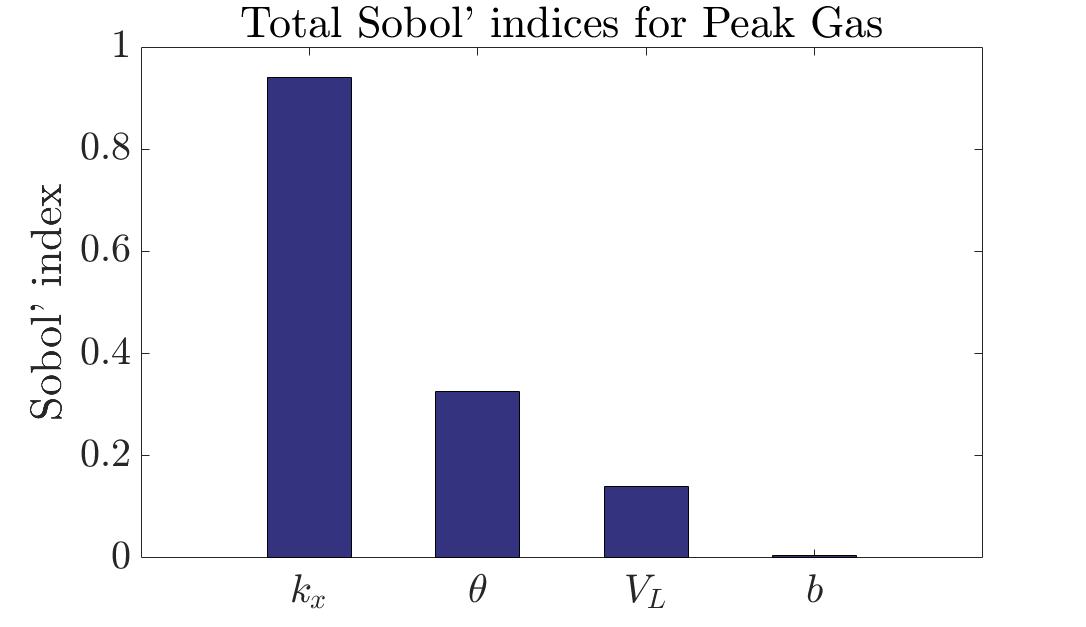}
\end{tabular}}
\end{center}
\caption{Bar charts of the estimated Sobol' indices.}
\label{fig:SA-bar-plots}
\end{figure}

As detailed in Section \ref{subsubsec:GSA} the surrogate PC models can be used to perform  global sensitivity analysis for the commercial model $\cM$. Using the Gauss-Legendre PCE with $p=6$ we estimate the Sobol' indices for both single inputs and pairs of inputs as well as the total Sobol' indices. The results are shown in Tables \ref{tab:sobol_singletons},  \ref{tab:sobol_pairs}  and \ref{tab:TEI} and the main effect Sobol' indices are illustrated in the bar charts of Figure \ref{fig:SA-bar-plots}.
Sudret \cite{Sudret2008} shows that for several analytic models even low order PCEs appear more reliable than Monte-Carlo based estimates; indicating that a comparison of Monte-Carlo estimates for the Sobol' indices of the non-analytic models against the estimates given by the PCEs may be an inappropriate measure. Thus, our discussion focusses on the Sobol' indices for the PCE surrogate models rather than the, unknown, Sobol' indices for the original models.

\begin{table}
\begin{center}
\resizebox{\linewidth}{!}{\begin{tabular}{lcrr}
\hline
\textbf{Variable}&\textbf{Sobol' Index} & \textbf{Total Gas} & \textbf{Peak Gas}\\
\hline
Fracture Permeability & $S_{(k_x)}$ & 0.1104156 & 0.7037193\\
Fracture Porosity & $S_{(\phi)}$ & 0.0056838 & 0.1583796\\
Langmuir Volume & $S_{(V_L)}$ & 0.8069139 & 0.0632314\\
Reciprocal of the Langmuir Pressure & $S_{(b)}$ & 0.0379921 & 0.0008090\\
\hline	
\end{tabular}}
\caption{Main effect Sobol' indices.}\label{tab:sobol_singletons}
\end{center}
\end{table}

\begin{table}
\begin{center}
\begin{tabular}{lcrr}
\hline
\textbf{Variables} & \textbf{Sobol' Index} & \textbf{Total Gas} & \textbf{Peak Gas}\\
\hline
Permeability and Porosity & $S_{(k_x,\phi)}$ & 0.0019166 & 0.0511791\\
Permeability and Volume & $S_{(k_x,V_L)}$ & 0.0294677 & 0.0209954\\
Permeability and Pressure & $S_{(k_x,b)}$ & 0.0027214 & 0.0002741\\
Porosity and Volume & $S_{(\phi,V_L)}$ & 0.0010781 & 0.0027350\\
Porosity and Pressure & $S_{(\phi,b)}$ & 0.0000810 & 0.0000342\\
Volume and Pressure & $S_{(V_L,b)}$ & 0.0052239 & 0.0000255\\
\hline	
\end{tabular}\\
\caption{Sobol' indices for pairwise interactions.}\label{tab:sobol_pairs}
\end{center}
\end{table}

\begin{table}
\begin{center}
\resizebox{\linewidth}{!}{\begin{tabular}{lcrr}
\hline
\textbf{Variable} & \textbf{Total Sobol' Index} & \textbf{Total Gas} & \textbf{Peak Gas}\\
\hline
Fracture Permeability & $T_{(k_x)}$ & 0.2162202 & 0.9411912\\
Fracture Porosity & $T_{(\phi)}$ & 0.0183316 & 0.3250164\\
Langmuir Volume & $T_{(V_L)}$ & 0.9607622 & 0.1397639\\
Reciprocal of the Langmuir Pressure & $T_{(b)}$ & 0.0679631 & 0.0043490\\
\hline	
\end{tabular}}
\caption{Total Sobol' indices.}\label{tab:TEI}
\end{center}
\end{table}

\begin{figure}
\begin{center}

%\resizebox{\linewidth}{!}{\begin{tabular}{ccc}
\includegraphics[width = 10cm]{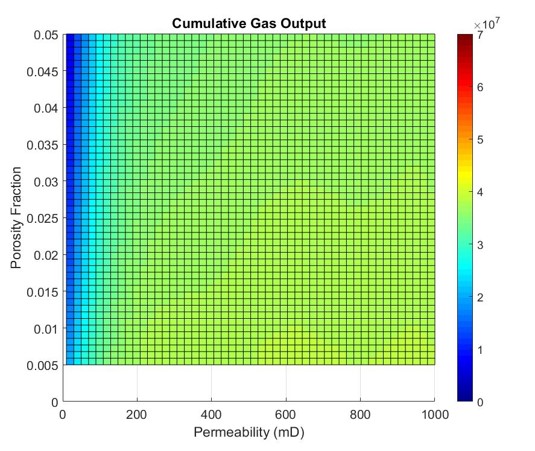} %& $\;$ & 
\includegraphics[width = 10cm]{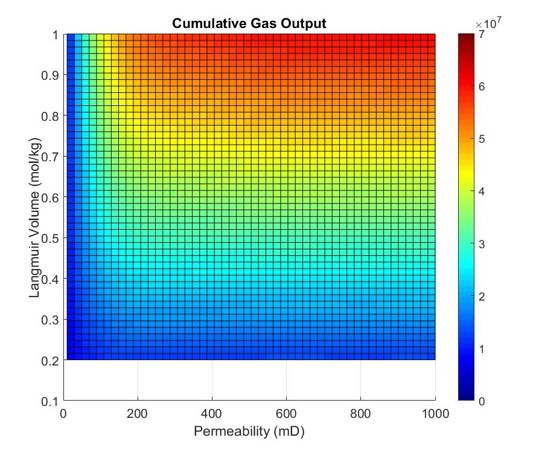}%\\
%\end{tabular}}
\end{center}
\caption{Slices of the PCE response surface generated for cumulative gas output.}
\label{fig:slices}
\end{figure}

Note that for cumulative gas extraction the variable with the most significant impact on the variance of the model is, 
as to be expected, 
the Langmuir volume. This is illustrated in Figure \ref{fig:slices} by taking two slices of the five dimensional response surface (four dimensions for the input parameters, one dimension for the response) generated by the PCE; the first has fixed Langmuir volume and pressure and shows how the cumulative gas varies with respect to fracture porosity and fracture permeability, and the second has fixed fracture porosity and Langmuir pressure and shows how the cumulative gas varies with respect to fracture permeability and Langmuir volume.  In the case of peak gas the dominant variable in terms of contribution to the overall variance is the fracture porosity.

The Sobol' indices indicate that variance in the (inverse) Langmuir pressure has, unexpectedly, low impact on the variance of the models. This is particularly surprising in the case of cumulative gas, where high $P_L$ is expected to strongly influence increased cumulative gas production in active wells \cite{Bahrami2015}. This discrepancy may be due to the fact that higher $P_L$ leads to faster recovery of gas per unit drop in pressure, while the influence of Langmuir volume on cumulative gas production affects long term extraction that may not be economically viable.  Similar reasoning can account for the discrepancy between the impact of the variance in permeability on cumulative and peak gas compared to the belief that, economically, permeability is the most important indicator of the viability of CSG production \cite{Bahrami2015}.

\section{Conclusions}
In this paper we have constructed Polynomial Chaos Expansions to act as surrogate models to a CMG commercial solver used to estimate peak and cumulative gas extraction from a coal seam gas well.

Each PCE propagates uncertainty in four input variables and builds a surrogate five dimensional response surface. The polynomial expansion delivers fast evaluations across the entire parameter space. For instance a PCE completes 3000 evaluations in under a second (on an Intel(R) Core(TM) i7-4770 CPU at 3.40GHz) whereas the commercial solver requires about five minutes per single evaluation. Being able to cheaply evaluate many points across the parameter space provides uncertainty quantification through the generation of summary statistics and empirical probability and cumulative density functions.
In addition analytic global sensitivity analyses can be formulated from PCEs at negligible additional computational cost.

The success of this process is demonstrated by the low discrepancy between the surrogates and original models. The constructed PCEs of order six approximate the commercial solver peak gas extraction rate with relative root mean square error of less than $0.01$ and the cumulative gas extraction with a relative root mean square error less than $0.06$.

In terms of the given application the PCE surrogate models provide good estimators for mean values and standard deviation of both the peak and cumulative gas productions predicted by the CMG model. Moreover, the cumulative distribution functions and the probability density functions for the CMG models can be approximated with a high degree of confidence by using the PCE prediction. Furthermore the PCE surrogate models allowed for the generation of Sobol' indices that indicate the most significant impact on the variance of the cumulative gas model is, as to be expected, the Langmuir volume and in the case of peak gas the dominant contributing input variable is the fracture porosity.

\section*{Acknowledgements}

The authors would like to thank the anonymous referees for their helpful comments and suggestions.

\bibliography{CMG_paper_preprint}

\end{document}